\documentclass[12pt,leqno]{article}
\usepackage{amsmath,amssymb,MnSymbol,calrsfs,url}
	
\DeclareMathOperator{\gph}{graph}
\DeclareMathOperator{\spt}{spt}

\DeclareMathOperator{\dist}{dist}
\DeclareMathOperator{\Clos}{Clos}
\DeclareMathOperator{\Bdry}{Bdry}

\def\res{\hbox{ {\vrule height .22cm}{\leaders\hrule\hskip.2cm} }}

\newcommand{\B}{\mathbf{B}}
\newcommand{\C}{\mathbf{C}}
\newcommand{\e}{\mathbf{e}}
\newcommand{\E}{\mathbf{E}}
\newcommand{\I}{\mathbf{I}}
\newcommand{\bL}{\mathbf{L}}
\newcommand{\M}{\mathbf{M}}
\newcommand{\N}{\mathbf{N}}
\newcommand{\p}{\mathbf{p}}
\newcommand{\R}{\mathbf{R}}
\newcommand{\bS}{\mathbf{S}}
\newcommand{\T}{\mathbf{T}}
\newcommand{\U}{\mathbf{U}}
\newcommand{\V}{\mathbf{V}}
\newcommand{\W}{\mathbf{W}}

\newcommand{\NU}{\boldsymbol{\nu}}
\newcommand{\MU}{\boldsymbol{\mu}}
\newcommand{\ALPHA}{\boldsymbol{\alpha}}
\newcommand{\GAMMA}{\boldsymbol{\gamma}}
\newcommand{\DELTA}{\boldsymbol{\delta}}
\newcommand{\BETA}{\boldsymbol{\beta}}

\newcommand{\HH}{\mathcal{H}}
\newcommand{\LL}{\mathcal{L}}

\newcommand{\RR}{\mathcal{R}}
\newcommand{\TT}{\mathcal{T}}

\newcommand{\graph}[2]{\gph_{#1} #2}

\numberwithin{equation}{section}
\newtheorem{theorem}[equation]{Theorem}
\newtheorem{lemma}[equation]{Lemma}
\newtheorem{corollary}[equation]{Corollary}

\begin{document}
\begin{flushleft}

TITLE: Partial boundary regularity for co-dimension one area-minimizing currents at immersed $C^{1,\alpha}$ tangential boundary points.

\medskip

ABSTRACT: We give partial boundary regularity for co-dimension one absolutely area-minimizing currents at points where the boundary consists of a sum of $C^{1,\alpha}$ submanifolds, possibly with multiplicity, meeting tangentially, given that the current has a tangent cone supported in a hyperplane with constant orientation vector; this partial regularity is such that we can conclude the tangent cone is unique. The proof follows closely the boundary regularity result given by Hardt and Simon in \cite{HS79}.

\medskip


\medskip


\section{Introduction}

Through a careful modification of the work found in \cite{HS79}, we are able to give partial regularity for co-dimension one absolutely area-minimizing currents at points where the boundary is tangentially $C^{1,\alpha}$ immersed. Our main result, Theorem \ref{main}, can be heuristically stated as follows:

\medskip

{\bf Theorem \ref{main}} {\it Suppose $T$ is an $n$-dimensional absolutely area-minimizing integer rectifiable current in an open subset of $\R^{n+1}$ containing the origin, and that near the origin $\partial T$ consists of a sum of $C^{1,\alpha}$ submanifolds for some $\alpha \in (0,1]$, each possibly with multiplicity, meeting tangentially (with same orientation) at the origin. Suppose as well that $T$ has a tangent cone at the origin $Q = M \E^{n} \res \{ (y_{1},\ldots,y_{n}) : y_{n} > 0 \} + m \E^{n} \res \{ (y_{1},\ldots,y_{n}): y_{n} < 0 \}$ where $M \geq 2$ and $m \in \{ 1,\ldots,M-1 \}.$ Then near the origin, there is a region of the horizontal hyperplane $\R^{n} \times \{0\}$ such that the support of $T$ over this region is the graph of a $C^{1,\beta}$ function for $\beta = \frac{\alpha}{4n+6}.$}

\medskip

Furthermore, the region is such that we can conclude $Q$ is the unique tangent cone of $T$ at the origin. Here, $\E^{n}$ is the current associated to the hyperplane $\R^{n} \times \{0\}$ with usual orientation, see 4.1.7 of \cite{F69}. See 4.3.16 of \cite{F69} for the definition of a tangent cone of a current.

\medskip

Theorem \ref{main} is precisely a generalization of Corollory 9.3 of \cite{HS79}, after applying the Hopf-type boundary point lemma given by Lemma 7  of \cite{FG57}, also appearing in \cite{HS79} as Lemma 10.1. We can get full boundary regularity via \cite{W83} in the special case that $\partial T$ is supported on exactly one $C^{1,\alpha}$ submanifold (if for example $m = M-1$), letting in this case $m \in \{0,\ldots,M-1\}$ and $M \geq 1.$ By \cite{W83} and the fact that the tangent cone of $T$ at the origin is $Q$ as above, if $m=0$ then $T$ corresponds to a $C^{1,\alpha}$ hypersurface-with-boundary, and if $m \geq 1$ then the support of $T$ near the origin is a real analytic hypersurface, with $T$ having multiplicity $M,m$ on either side of $\partial T.$     

\subsection{An application of Theorem \ref{main}} 

We note an application of Theorem \ref{main}, which in fact motivated the present work. Recently in \cite{R13} the author introduced the $c$-isoperimetric mass of currents, which is defined for each $c>0$ by $$\M^{c}(T) = \M(T) + c \M (\partial T)^{\kappa}$$ whenever $T$ is an $n$-dimensional integer multiplicity rectifiable current in $\R^{n+k},$ $\M$ is the usual mass on currents, and $\kappa = \frac{n}{n-1}$ is the isoperimetric exponent. 

\medskip

This leads to define and study a minimization problem. Let $\Gamma$ be an $(n-1)$-dimensional integer rectifiable current in $\R^{n+k}$ with compact support and $\partial \Gamma = 0,$ which we refer to as the fixed boundary. Define $\I_{\Gamma}(\R^{n+k})$ to be the set of $n$-dimensional integer rectifiable currents $T$ with compact support so that $\partial T = \Gamma + \Sigma$ where $\Gamma$ and $\Sigma$ have disjoint supports. We then say $\T_{c} \in \I_{\Gamma}(\R^{n+k})$ is a solution to the $c$-Plateau problem with respect to fixed boundary $\Gamma$ if $\T_{c}$ minimizes $\M^{c}$ amongst all $T \in \I_{\Gamma}(\R^{n+k})$ (see Definition 3.3 of \cite{R13} with $U = \R^{n+k}$). For such $\T_{c},$ writing $\partial \T_{c} = \Gamma + \Sigma_{c}$ we refer to $\Sigma_{c}$ as the free boundary.

\medskip

Theorem 8.2 of \cite{R13} concludes there is no solution to the $c$-Plateau problem $\T_{c}$ with $\partial \T_{c} = \Gamma + \Sigma_{c}$ with nonzero free boundary $\Sigma_{c}$ a smooth embedded $(n-1)$-dimensional submanifold with parallel mean curvature, that is constant mean curvature in the sense of \cite{H73}, so that $\T_{c}$ near $\Sigma_{c}$ is a smooth submanifold-with-boundary. This can be used in Theorem 9.1 of \cite{R13} to show that in case the fixed boundary $\Gamma$ is one-dimensional in the plane, that is if $n=2,k=0,$ then free boundaries must always be empty. However, so-called non-trivial solutions in the limit can occur, as seen in Theorem 10.2 of \cite{R13} which shows that for small values of $c>0$ when $\Gamma$ is the square in the plane, the infimum of $\M^{c}$ is attained in the limit by a sequence of currents in $\I_{\Gamma}(\R^{2})$ which converge to a nonempty current not in $\I_{\Gamma}(\R^{2}).$

\medskip

The author conjectures that this holds generally in $n=2,k=1:$ if the fixed boundary $\Gamma$ is one-dimensional in $\R^{3},$ then for each $c>0$ either every solution to the $c$-Plateau problem $\T_{c}$ with fixed boundary $\Gamma$ has empty free boundary, so that $\partial \T_{c} = \Gamma,$ or the infimum value of $\M^{c}$ can only be attained in the limit by a sequence of currents in $\I_{\Gamma}(\R^{3}).$ Via preliminary sketches of arguments, the author strongly conjectures that one-dimensional free boundaries in space are always smooth, so that by Theorem 8.2 of \cite{R13} it only suffices to show they have parallel mean curvature. We claim this follows through a geometric analysis, by arguments similar to those used in \cite{R10} and \cite{R11} to study the two-valued minimal surface equation, a PDE first introduced in \cite{SW07}. The partial regularity result given by Theorem \ref{main} is essential in order to carry over this analysis. 

\subsection{Counterexamples and questions}

The examples of stable branched minimal immersions given by \cite{SW07} and \cite{R10} show the absolutely area-minimizing hypothesis cannot be relaxed to stability. Indeed, Theorem 1 of \cite{SW07} holds that if $u_{0}$ is a solution to the two-valued minimal surface equation (see the operator $\mathcal{M}_{0}$ at the start of \S 3 of \cite{R10}) over the punctured unit disk in $\R^{2}$ which can be extended continuously across the origin, then $$G = \{ (re^{i \theta},u_{0}(r^{1/2} e^{i \theta/2})): r \in (0,1), \theta \in \R \}$$ is a stable minimal immersion with $C^{1,\alpha}$ branch point at $(0,u_{0}(0)),$ for some $\alpha \in (0,1).$ \cite{SW07} and \cite{R10} show a large non-trivial class of such solutions exist. We can thus show there is a solution $u_{0}$ to the two-valued minimal surface equation which can be extended continuously across the origin, so that $\{ (re^{i \theta},u(r^{1/2}e^{i\theta/2})): r \in (0,1), \theta \in (0,3\pi) \}$ satisfies the assumptions of Theorem \ref{main} (with $M=2,m=1,$ and with the absolutely area-minimizing condition replaced by stability) but fails to satisfy the partial regularity conclusions given there.  

\medskip

Neither does Theorem \ref{main} hold in higher co-dimensions. A counterexample is given by considering the region $\{ (re^{i \theta}, r^{3/2} e^{\frac{3i \theta}{2}}) : r > 0, \ \theta \in [0,3 \pi] \}$ of the holomorphic variety $\{ (z,w) : z^{3} = w^{2} \} \subset \C \times \C \cong \R^{4},$ which is still calibrated and hence area-minimizing. The best general result in all co-dimensions is thus as in \cite{A75}, boundary regularity in case of currents with $C^{1,\alpha}$ embedded boundary at points of density near $1/2$; see Theorem 0.1 of \cite{DS02}, which concludes this in fact for almost minimizing currents of arbitrary co-dimension, or more generally, \cite{B10} which does this for stationary varifolds. Observe again, that the examples from \cite{SW07} and \cite{R10} show the density $=1/2$ assumption cannot be relaxed without the area-minimizing hypothesis.

\subsection{Summary}

The aim of this work is thus to extend Corollary 9.3 of \cite{HS79} to the conditions set forth by Theorem \ref{main}, in order to study the $c$-Plateau problem in space as introduced in \cite{R13}. We now discuss the organization of this work. 

\medskip

The proof of Theorem \ref{main} involves making small but ubiquitous changes to the proofs found in \cite{HS79}. This task is undertaken in \S 3, which we begin with a discussion in order to facilitate the reader in making those modifications to \cite{HS79}. Nevertheless, we strongly recommend that the reader be familiar with the arguments of \cite{HS79}. Early in \S 3 we include a list of notation used throughout.

\medskip

To modify \cite{HS79}, we must rely on the calculations established in the Appendix, which contain the deeper differences between the present setting and the proof of \cite{HS79}. The key identity, which demonstrates why we can modify \cite{HS79}, can be found in \eqref{appendixA6}. From this calculation one can conclude that if $T$ satisfies the assumptions of Theorem \ref{main}, then for $\p$ the projection onto the subspace $\R^{n} \times \{0\}$ we have that near the origin $\p_{\#} T - \E^{n}$ and $\E^{n}$ have additive masses.

\medskip

Before all this, we state in \S 2 the main result Theorem \ref{main}, giving exactly the assumptions necessary. Corollary \ref{uniquetangentcones} concludes uniqueness of tangent cones for $T$ satisfying the conditions of Theorem \ref{main}. We also remark in Theorem \ref{existencetangentcones}, using \cite{B76}, why currents $T$ as in Theorem \ref{main} must have tangent cones at such tangential boundary points.   

\section{Main Results}

For the definition of absolutely area-minimizing, consult 5.1.6 of \cite{F69}. We denote $\E^{n-1}, \E^{n}$ to be the currents associated respectively to $\R^{n-1} \times \{0\},\R^{n} \times \{0\}$ in $\R^{n+1}$ each with usual orientation, as in 4.1.7 \cite{F69}. Given $r>0,$ we define the map $\MU_{r}(x) = rx,$ and for a current $T$ we let $\MU_{r \#} T$ be the push-forward of $T$ by $\MU_{r}.$ Let also $\Clos A$ denote the closure of $A \subset \R^{n+1}.$ We now state our main result. 

\begin{theorem} \label{main} Suppose $\alpha \in (0,1]$ and $T$ is an $n$-dimensional absolutely area-minimizing locally rectifiable integer multiplicity current in $\R^{n+1} \cap \{ x:|x| < 3\}.$ We also suppose $T$ satisfies the hypothesis:

\begin{enumerate} \item[($\ast$)] $\begin{aligned} \partial T \res \{ (x_{1},\ldots,x_{n+1}): & |(x_{1},\ldots,x_{n-1})| < 2, \ |x_{n}| < 2 \} \\ & = (-1)^{n} \sum_{\ell=1}^{\N} m_{\ell} \Phi_{T,\ell \#}(\E^{n-1} \res \{ z: |z| < 2 \}), \end{aligned}$ \\ where $m_{\ell}$ are positive integers, and for each $\ell \in \{1,\ldots,\N\}$ $$\Phi_{T,\ell}(z_{1},\ldots,z_{n-1}) = (z_{1},\ldots,z_{n-1},\varphi_{T,\ell}(z_{1},\ldots,z_{n-1}),\psi_{T,\ell}(z_{1},\ldots,z_{n-1}))$$ where $\varphi_{T,\ell},\psi_{T,\ell} \in C^{1,\alpha}(\R^{n-1} \cap \{ z:|z|<2 \})$ with $\varphi_{T,\ell}(0) = 0 = \psi_{T,\ell}(0), \ D \varphi_{T,\ell}(0) = 0 = D \psi_{T,\ell}(0).$
 \item[($\ast \ast$)] $T$ has a tangent cone $$\Big[ M \E^{n} \res \{ (y_{1},\ldots,y_{n}): y_{n} > 0 \} + m \E^{n} \res \{ (y_{1},\ldots,y_{n}): y_{n} < 0 \} \Big] \times \delta_{0}$$ at the origin, where $M \in \{ 2,3,\ldots \}$ and $m \in \{ 1,\ldots,M-1 \}.$
\end{enumerate}

Then there is a $\delta = \delta(n,m,M,\alpha) \in (0,1)$ sufficiently small, so that letting $$\begin{aligned} \tilde{V} & = \{ y=(y_{1},\ldots y_{n}) : y_{n} > |y|^{1+\beta}, \ |y| < \delta \} \\ \tilde{W} & = \{ y=(y_{1},\ldots y_{n}) : y_{n} < - |y|^{1+\beta}, \ |y| < \delta \} \end{aligned}$$ for $\beta = \alpha/(4n+6),$ then for $\rho>0$ sufficiently small depending on $T$ we have $$\begin{aligned} \p^{-1}(\tilde{V}) \cap \spt \MU_{1/\rho \#} T & = \graph{\tilde{V}}{\tilde{v}} \\ \p^{-1}(\tilde{W}) \cap \spt \MU_{1/\rho \#} T & = \graph{\tilde{W}}{\tilde{w}} \end{aligned}$$ for some $\tilde{v} \in C^{1,\beta}(\Clos \tilde{V}), \tilde{w} \in C^{1,\beta}(\Clos \tilde{W})$ such that $\tilde{v}|\tilde{V}, \tilde{w}|\tilde{W}$ satisfy the minimal surface equation and $D \tilde{v}(0) = 0 = D \tilde{w}(0).$ Furthermore, we have $$\sup_{y \in \tilde{V}} \frac{|D^{2}\tilde{v}(y)|}{|y|^{\beta-1}} + \sup_{y \in \tilde{W}} \frac{|D^{2}\tilde{w}(y)|}{|y|^{\beta-1}} \leq c$$ for some $c=c(n,m,M) \in (0,\infty).$ \end{theorem}

Note that $M-m = \sum_{\ell=1}^{\N} m_{\ell}.$ As noted in the introduction, the case $m=M-1$ is just Corollary 9.3 of \cite{HS79}, together with Lemma 10.1 of \cite{HS79}. Also, if $\N=1,$ that is when $\partial T$ is a $C^{1,\alpha}$ $(n-1)$-dimensional submanifold with multiplicity, then Theorem \ref{main} follows in this case by the higher multiplicity boundary regularity given by \cite{W83}. Nonetheless, the proof we give below will cover all cases. The following corollary immediately follows.

\begin{corollary} \label{uniquetangentcones} If $T$ is as in Theorem \ref{main}, then $T$ has unique tangent cone $$M \E^{n} \res \{ (y_{1},\ldots,y_{n}): y_{n} > 0 \} + m \E^{n} \res \{ (y_{1},\ldots,y_{n}): y_{n} < 0 \}$$ at the origin. \end{corollary}

Before proceeding, we prove a lemma showing the existence of tangent cones to start.

\begin{theorem} \label{existencetangentcones} Suppose $\alpha \in (0,1]$ and $T$ is an $n$-dimensional absolutely area-minimizing locally rectifiable integer multiplicity current in $\R^{n+1} \cap \{ x: |x| < 3 \}$ satisfying hypothesis $(\ast).$ Then $T$ has an oriented tangent cone at the origin, and every oriented tangent cone of $T$ at the origin is absolutely area minimizing with density at the origin equal to the density of $T$ at the origin. \end{theorem}

{\it Proof.} By Theorems 3.6,3.3 of \cite{B76} we only need to check the finiteness of $$\nu_{1}^{\partial T}(x) = \int_{\R^{n+1} \cap \{ x: |x| < 1 \}} \frac{|\vec{\partial T} \wedge (y-x)|}{|y-x|^{n}} \ d \| \partial T \| (y).$$ This however follows by hypothesis $(\ast).$

\medskip

Observe of course that $T$ may satisfy hypothesis $(\ast)$ but not ($\ast \ast$), if for example $T$ is a union of half-planes in space, appropriately oriented.

\section{Proof of Theorem \ref{main}}

We must follow closely the arguments of \cite{HS79}, up to Corollary 9.3 found therein. Using the calculations of the Appendix, we must make clear to the reader where and what changes must be made to fit the setting given by Theorem \ref{main}. In \S 3.1, we modify the necessary results of \cite{HS79}, and in \S 3.2 we finish the proof of Theorem \ref{main}.

\subsection{Modifying the results of \cite{HS79}}

There are twelve sections in \cite{HS79}, each devoted to a theoretical step. Each section of \cite{HS79} is further divided into subsections, given either by closely related computations, lemma, or theorem. We refer to the sections and subsections of \cite{HS79} with the marker ``HS"; the sections of \cite{HS79} we must discuss in detail are HS1-HS9. However, we adopt the rule that when a serious difference must be made to a subsection of \cite{HS79}, then that subsection is denoted by the prefix ``R"; for example, we must make a notable change to the lemma found in HS3.2, and so we refer instead to R3.2 in what follows. We avoid restating whole lemmas or theorems from \cite{HS79} when the only difference is due to notation (owing to the present setting involving different multiplicities).

\medskip

Constants $c_{1},\ldots,c_{49}$ are introduced in \cite{HS79}. Our underlying goal is to show there are analogous constants here depending on $m,M,n$ (in fact, we will only need to discuss constants $c_{1},\ldots,c_{46}$). This shall be crucial in applying many of the iterative arguments found in \cite{HS79} to the present setting. We include minor, although clarifying, corrections to \cite{HS79}. We start each section of \cite{HS79} with a general description for convenience to the reader. Nonetheless, we shall be as succinct as possible.

\medskip

\begin{center}{\bf HS1. Notation and preliminaries}\end{center}

\medskip

This section establishes the basic notation used throughout, which we use unchanged, the only difference being in R1.6. We fix in this case $$M \in \{ 2,3,\ldots \}, \ m \in \{ 1,\ldots,M-1 \}, \ n \in \{ 2,3,\ldots \}, \ \alpha \in (0,1].$$ Our goal is to find constants $c_{1},\ldots,c_{49},$ playing analogous roles to those found in \cite{HS79}, which depend only on $m,M,n.$  

\medskip

HS1.1. {\it Standard notation.} We follow the standard notation found in 668-671 of \cite{F69}. Notably is $$\Clos A, \ \Bdry A, \ A \sim B, \ \R^{n}, \ \LL^{n}, \ \RR_{n}(\R^{n+1}),$$ $$\E^{k}, \ [a,b], \ \DELTA_{a}, \ \partial T, \ f_{\#} T, \ T \res A, \ \M(T), \ \spt T, \ \|T\|, \ \vec{T}, \ \Theta^{n}(\|T\|,a).$$

\medskip

HS1.2. {\it Special notation associated with} $\R^{n}$. We use the same notation, which is as follows: $$\begin{aligned} \U^{n}(y,r) & = \R^{n} \cap \{ z: |y-z| < r \} \text{ and } \\ \B^{n}(y,r) & = \R^{n} \cap \{ z: |y-z| \leq r \} \text{ for } y \in \R^{n} \text{ and } 0 < r < \infty, \\ \ALPHA(n) & = \LL^{n}[\U^{n}(0,1)], \\ \bL & = \U^{n}(0,1) \cap \{ (y_{1},\ldots,y_{n}):y_{n}=0 \} \\ \V & = \U^{n}(0,1) \cap \{ (y_{1},\ldots,y_{n}):y_{n}>0 \} \\ \W & = \U^{n}(0,1) \cap \{ (y_{1},\ldots,y_{n}):y_{n}<0 \} \\ \V_{\sigma} & = \V \cap \{ y: \dist(y,\Bdry \V) > \sigma \} \text{ and } \\ \W_{\sigma} & = \W \cap \{ y: \dist(y, \Bdry \W) > \sigma \} \text{ for } 0<\sigma<1. \end{aligned}$$
 
\medskip

HS1.3. {\it Special notation associated with} $\R^{n+1}.$ We use the following: $$\begin{aligned} \U_{r} & = \R^{n+1} \cap \{ x:|x| < r \}, \\ \B_{r} & = \R^{n+1} \cap \{ x:|x| \leq r \}, \\ \C_{r} & = \R^{n+1} \cap \{ x:|\p(x)| \leq r \} \text{ where } 0<r<\infty \text{ and } \\ & \p : \R^{n+1} \rightarrow \R^{n}, \ \p(x_{1},\ldots,x_{n+1}) = (x_{1},\ldots,x_{n}). \end{aligned}$$ $$\begin{aligned} \e_{1} = (1,0,\ldots,0) \in \R^{n+1}, \e_{2} & = (0,1,0,\ldots,0) \in \R^{n+1},\ldots, \\ \e_{n+1} & = (0,\ldots,0,1) \in \R^{n+1}, \end{aligned}$$ $$X_{k}:\R^{n+1} \rightarrow \R, \ X_{k}(x)=x_{k}$$ $$X = (X_{1},\ldots,X_{n+1}) \text{ (the identity on $\R^{n+1}$) }, \ |X| = (\sum_{k=1}^{n+1} X_{k}^{2})^{1/2},$$ $$\begin{aligned} & \MU_{r}: \R^{n+1} \rightarrow \R^{n+1}, \ \MU_{r}(x) = rx, \\ & \BETA_{r}: \R^{n+1} \rightarrow \R^{n+1}, \ \BETA_{r}(x) = (x_{1},\ldots,x_{n},rx_{n+1}), \\ & \GAMMA_{\omega}: \R^{n+1} \rightarrow \R^{n+1}, \\ \GAMMA_{\omega}(x) = (x_{1},& \ldots,x_{n-1},x_{n} \cos \omega - x_{n+1} \sin \omega,x_{n} \sin \omega+x_{n+1} \cos \omega),\end{aligned}$$ for $x=(x_{1},\ldots,x_{n+1}) \in \R^{n+1},$ $k \in \{1,\ldots,n+1\},$ $0<r<\infty,$ and $\omega \in \R.$

\medskip

HS1.4. {\it Special notation associated with} $T \in \RR_{n}(\R^{n+1}).$ For $T \in \RR_{n}(\R^{n+1}),$ let $$\NU^{T} = (\NU^{T}_{1},\ldots,\NU^{T}_{n+1})$$ be the {\it unit normal vectorfield} associated with $T,$ where we define $\NU^{T}_{k}$ by $\vec{T} = \sum_{k=1}^{n+1} \NU^{T}_{k} (-1)^{k} \e_{1} \wedge \ldots \wedge \e_{k-1} \wedge \e_{k+1} \wedge \ldots \wedge \e_{n+1}$ (see also \cite{HS79}), and $$\delta^{T} = (\delta^{T}_{1},\ldots,\delta^{T}_{n+1})$$ be the {\it tangential gradient operator} associated with $T$ so that, for $\|T\|$ almost all $x \in \R^{n+1},$ $$\delta^{T} f(x) = Df(x) - [\NU^{T}(x) \cdot Df(x)] \NU^{T}(x) \text{ for } f \in C^{1}(\R^{n+1})$$ (where $Df(x) = (D_{1}f(x),\ldots,D_{n+1}f(x)) \in \R^{n+1}$) is the {\it tangent gradient} of $f,$ $$\delta^{T} \cdot g(x) = \sum_{k=1}^{n+1} \delta_{k}^{T} g_{k}(x) \text{ for } g_{1},\ldots,g_{n+1} \in C^{1}(\R^{n+1}) \text{ and } g=(g_{1},\ldots,g_{n+1})$$ is the {\it tangential divergence of} $g.$

\medskip

Of central importance, as in \cite{HS79}, are for $0<r<\infty$ $$\E_{\bS}(T,r) = r^{-n} \M(T\res\B_{r}) - \ALPHA(n) \Theta^{n}(\|T\|,0)$$ the {\it spherical excess} (whenever $\Theta^{n}(\|T\|,0)$ exists), and $$\E_{\C}(T,r) = r^{-n} \M(T\res\C_{r}) - r^{-n}\M\p_{\#}(T\res\C_{r})$$ the {\it cylindrical excess}.

\medskip

Observe that HS1.4(1) should in fact be, for $0<r<s<\infty$ $$\begin{aligned} r^{n} \E_{\C}(T,r) & = \int_{\C_{r}} [1-|\NU^{T} \cdot \e_{n+1}|] d \| T \| \\ & \leq \int_{\C_{s}} [1-|\NU^{T} \cdot \e_{n+1}|] \ d \| T \| = s^{n} \E_{\C}(T,s). \end{aligned}$$

\medskip

R1.5. {\it Tangent cones at the boundary.} This section is supplanted exactly by the hypothesis ($\ast \ast$) of Theorem \ref{main}.

\medskip

R1.6. {\it The family} $\TT$. We must account for the different setting given by Theorem \ref{main}. We define the family $\TT$ in this case as follows (having fixed $m,M$ as above).

\medskip

{\it Let $\TT$ denote the collection of all absolutely area minimizing $T \in \RR_{n}(\R^{n+1})$ such that} $$\begin{aligned} \spt T & \subset B_{3} \\ \M(T) & \leq 3^{n} [1+ M \ALPHA (n)] \\ \Theta^{n}(\|T\|,0) & = \frac{M+m}{2}. \end{aligned}$$ {\it Also, $T$ satisfies hypothesis $(\ast)$ (see the statement of Theorem \ref{main}) with} $$\begin{aligned} \kappa_{T} = & 2 \alpha^{-1} \max_{\ell=1,\ldots,N} \sup_{z \neq w} |z-w|^{-\alpha} [|D \varphi_{T,\ell}(z)-D \varphi_{T,\ell}(w)|^{2} \\ & +|D \psi_{T,\ell}(z)-D \psi_{T,\ell}(w)|^{2}]^{1/2} \leq 1, \end{aligned}$$ {\it and if we define $\varphi^{max}_{T},\varphi^{min}_{T}: \R^{n-1} \cap \{ z : |z| \leq 2 \} \rightarrow \R$ by} $$\varphi^{max}_{T}(z) = \max_{\ell = 1,\ldots,\N} \varphi_{T,\ell}(z), \hspace{.25in} \ \varphi^{min}_{T}(z) = \min_{\ell = 1,\ldots,\N} \varphi_{T,\ell}(z),$$ {\it then} $$\begin{aligned} \p_{\#} (T \res \C_{2}) \res & \{ (y_{1},\ldots,y_{n}): y_{n} > \varphi_{T}^{max}(y_{1},\ldots,y_{n-1}) \text{ or } y_{n} < \varphi_{T}^{min}(y_{1},\ldots,y_{n-1}) \} \\ & = M [\E^{n} \res \U^{n}(0,2) \cap \{ (y_{1},\ldots,y_{n}): y_{n} > \varphi^{max}_{T}(y_{1},\ldots,y_{n-1}) \}] \\ & \ + m [\E^{n} \res \U^{n}(0,2) \cap \{ (y_{1},\ldots,y_{n}): y_{n} < \varphi^{min}_{T}(y_{1},\ldots,y_{n-1}) \}]. \end{aligned}$$ We see that we must modify HS1.6(1), using \eqref{appendixA5},\eqref{appendixA4} we get for $r \in (0,2]$ \begin{equation} \tag*{R1.6(1)} \begin{aligned} \E_{\bS}(T,r) & \leq r^{-n} \| T \| \C_{r} - \left( \frac{M+m}{2} \right) \ALPHA (n) \\ & \leq \E_{\C}(T,r) +r^{-n} \M \p_{\#} (T \res \C_{r}) - \left( \frac{M+m}{2} \right) \ALPHA (n) \\ & \leq \E_{\C}(T,r) + (M-m) \ALPHA (n-1) r^{\alpha} \kappa_{T}  \end{aligned} \end{equation}

\medskip

\begin{center}{\bf HS2. First variation and monotonicity}\end{center}

\medskip

Monotonicity formulas are computed in this section, via the first variation, which hold here without major change. We introduce in this section $c_{1},\ldots,c_{5}.$

\medskip

HS2.1. {\it First variation.} The first variation formula of course holds here as stated.

\medskip

HS2.2. {\it Monotonicity estimates} The formulas HS2.2(1)-(5) holds with $c_{1},\ldots,c_{4}$ depending on $m,M,n.$ This follows via two observations.

\medskip

First, is that by $(\ast)$ we have $$2^{n} \int |g \wedge \vec{\partial T}| \ d \| \partial T \| = 2^{n} \sum_{\ell=1}^{\N} m_{\ell} \int_{\Phi_{T,\ell}(\R^{n-1} \cap \{z: |z| < 2\})} |x \wedge \vec{\partial T}| \ d \HH^{n-1}.$$ Furthermore, HS2.2(6) holds with $c_{5}$ as in \cite{HS79} for $\HH^{n-1}$ almost-every $x \in \Phi_{T,\ell}(\R^{n-1} \cap \{z: |z| < 2\}),$ for each $\ell = 1,\ldots,\N.$  

\medskip

Second, define $L:\R^{n+1} \rightarrow \R^{n}$ as in \cite{HS79} by $$L(x_{1},\ldots,x_{n+1}) = ((x_{1},\ldots,x_{n-1}),|(x_{n},x_{n+1})|) \text{ for } (x_{1},\ldots,x_{n+1}) \in \R^{n+1}.$$ We still have for $0 < r < 2$ that $$r^{-n} \| L_{\#} T \| \B_{r} \geq r^{-n} \M \E^{n} \res \{ y = (y_{1},\ldots,y_{n}): 2|(y_{1},\ldots,y_{n-1})| \leq |y| \leq r, \ y_{n}>0 \}.$$ We can then follow the arguments used to show HS2.2(1)-(5). However, we must allow $c_{1},\ldots,c_{4}$ to depend on $m,M$ on account of the estimate on $2^{n} \int |g \wedge \vec{\partial T}| \ d \| T \|$ using $(\ast)$ in showing HS2.2(4).

\medskip

HS2.3. {\it Remark.} With $T \in \TT$ and $3 \leq r < \infty,$ using HS2.2(4) and R1.6(1) $$\M[(\MU_{r \#}T) \res B_{3}] < 3^{n} \left[ E_{\C}(T,1) + \left( \frac{M+m}{2} \right) \ALPHA (n) + ((M-m) \ALPHA (n-1) + c_{4}) \kappa_{T} \right].$$ We can then choose $c_{4}$ depending on $m,M,n$ so that $\E_{\C}(T,1) + \kappa_{T} \leq (1+c_{4})^{-1}$ implies $\M[(\MU_{r \#}T) \res B_{3}] < 3^{n} [1+ M \ALPHA (n)].$ From this HS2.3(1)(2) hold with no change.

\medskip

\begin{center}{\bf HS3. An area comparison lemma}\end{center}

\medskip

The results of this section shall be used in the next to conclude preliminary bounds on the excess. Although we must make a serious change to HS3.2, our version is sufficient. Observe that the last constant introduced in HS2 was $c_{5},$ in HS2.2, and the first introduced in HS3 is $c_{7},$ in R3.2. We introduce $c_{7},c_{8},c_{9},$ and unlike in \cite{HS79} have no need for $c_{10},c_{11}$ in R3.2.

\medskip

HS3.1. {\it Remark.} A general fact about exterior algebras is stated.

\medskip

R3.2. LEMMA. The conclusions of this lemma are different, and so we state the new version.

\medskip

{\it There are constants $c_{7},c_{8}$ depending on $m,M,n$ such that if $T \in \TT,$ $\rho > 0,$ $0 < \tau < 1,$ $A = \p^{-1}[\p(A)]$ is a Borel subset of $\C_{1},$ $A_{\tau} = \{ x:\dist (x,A) < \tau\},$ $\mu:\R^{n} \rightarrow \R$ is a $C^{1}$ function, $\sup_{\p(A)}|D \mu| \leq \rho/\tau$ and} $$F(x,_{1},\ldots,x_{n+1}) = (x_{1},\ldots,x_{n},\mu(x_{1},\ldots,x_{n}) x_{n+1})$$ {\it for $(x_{1},\ldots,x_{n+1}) \in \R^{n+1}$, then} $$\M F_{\#}(T \res A) - \M (T \res A) \leq c_{7} \kappa_{T} + c_{8} \tau^{-2} (1+\rho^{2}) \int_{A_{\tau}} X_{n+1}^{2} \ d \|T\|.$$

\medskip

{\it Proof.} Most of the proof is the same. The same calculation shows \begin{equation} \tag*{HS3.2(1)} \begin{aligned} \M & F_{\#}(T \res A) - \M(T \res A) \\ & \leq 2 \int_{A} [1- (\NU_{n+1}^{T})^{2}] \ d \| T\| + 2 \int_{A} X_{n+1}^{2} |D \mu|^{2} \ d \| T \|. \end{aligned} \end{equation} The second term is $\leq \frac{2 \rho^{2}}{\tau^{2}} \int_{A_{\tau}} X_{n+1}^{2} \ d \| T \|$ as is necessary.

\medskip

Taking the vector field $g = X_{n+1} \lambda^{2} \e_{n+1}$ into the first variation formula HS2.1, with $\lambda : \R^{n+1} \rightarrow [0,1]$ a $C^{1}$ function with $\spt \lambda \subset A_{\tau},$ $\lambda|_{A} \equiv 1,$ and $\sup |D \lambda| \leq c_{9}/\tau,$ yields $$\begin{aligned} \int [1 & -(\NU_{n+1}^{T})^{2}] \lambda^{2} \ d \|T\| \\ & = \int -2 \lambda X_{n+1} \delta^{T} \lambda \cdot \e_{n+1} \ d \|T\| + \int X_{n+1} \lambda^{2} \e_{n+1} \cdot \vec{\partial T} \ d \| \partial T \| \\ & \leq \frac{1}{2} \int (1-(\NU_{n+1}^{T})^{2}) \ d \|T\| + 2 \int X_{n+1}^{2} |D \lambda|^{2} \ d \|T\| + \int |X_{n+1}| \lambda^{2} \ d \| \partial T \|. \end{aligned}$$ We thus have $$\int [1-(\NU_{n+1}^{T})^{2}] \lambda^{2} \ d \|T\| \leq \frac{4 c_{9}^{2}}{\tau^{2}} \int_{A_{\tau}} X_{n+1}^{2} \ d \|T\| + 2 \int |X_{n+1}| \lambda^{2} \ d \| \partial T \|.$$ We can also compute $$\begin{aligned} \int |X_{n+1}| \lambda^{2} \ d \| \partial T \| & \leq \int_{\C_{1+\tau}} |X_{n+1}| \ d \| \partial T \| \\ & \leq \left(\frac{\alpha}{2}\right) \kappa_{T} (1+\tau)^{1+\alpha} \| \partial T \| \C_{1+\tau} \\ & \leq 2^{n-1+\alpha} \sqrt{3} (M-m) \ALPHA(n-1) \kappa_{T} \end{aligned}$$ using $\alpha,\kappa_{T},\tau \leq 1$ and $$\| \partial T \| \C_{1+\tau} \leq (M-m) \ALPHA(n-1) \left( 1+ \frac{\alpha^{2} \kappa_{T}^{2}}{4} (1+\tau)^{2 \alpha} + \frac{\alpha^{4} \kappa_{T}^{4}}{16} (1+\tau)^{4\alpha} \right)^{1/2} (1+\tau)^{n-1}.$$ We conclude R3.2 with $c_{7}$ depending on $m,M,n,$ and $c_{8}$ actually just depending on $c_{9}.$
 
\medskip

\begin{center}{\bf HS4. Some preliminary bounds on excess}\end{center}

\medskip

This section compares the cylindrical excess to the height excess, using subharmonicity while referring to either 7.5(6) of \cite{A72} or 3.4 of \cite{MS73}. The proofs and results are the same, although make a slight clarification to the proof of HS4.1. In this section we introduce $c_{12},\ldots,c_{16}.$

\medskip

HS4.1. LEMMA. The result is the same, with HS4.1(1)(2) holding exactly, although $c_{12},c_{13},c_{14},c_{15}$ depend on $m,M,n.$

\medskip

{\it Proof.} The second inequality in HS4.1(1) follows immediately from R1.6 with $c_{13}$ now depending on $M,n.$ We make a small clarification to the proof of the first inequality in HS4.1(1).

\medskip

Observe that if $\kappa_{T} > 3^{n}[1+M\ALPHA(n)] \sigma^{2},$ then by R1.6 $$c_{12}^{-1} \sigma^{2} \E_{\C}(T,1) - \kappa_{T} < 0$$ so long as we choose $c_{12} \geq 1.$ 

\medskip

We can assume $\kappa_{T} \leq 3^{n} [1+M \ALPHA(n)] \sigma^{2}.$ Letting $\tau = \sigma/2,$ $A = \C_{1+\tau} \sim \C_{1},$ and $\mu,$ $F,$ $h,$ and $R_{T}$ be as in \eqref{appendixA8} of Lemma \ref{appendixlemma3}, then as in HS3.2 $$\begin{aligned} 
\E_{\C}(T,1) & \leq \M F_{\#}(T \res A) - \M(T \res A) + \M (R_{T}) - \M F_{\#} (T \res \C_{1}) \\ & \leq \M F_{\#}(T \res A) - \M(T \res A) + \M (R_{T}) \end{aligned}$$ Using \eqref{appendixA8} (since $\kappa_{T} \leq 4 \cdot 3^{n} [1+M \ALPHA(n)] \tau^{2}$) and R3.2 with $\rho = 3$ gives $$\begin{aligned} \E_{\C}(T,1) \leq & c_{7} \kappa_{T} + 10c_{8} \tau^{-2} \int_{\C_{1+\sigma}} X_{n+1}^{2} \ d \|T\| \\ & + \frac{1}{2} \left[ \frac{\sqrt{21}}{4}+ 2^{\frac{9n-5}{2}} 3^{n^{2}-\frac{1}{2}} \right] (M-m) \ALPHA(n-1)[1+M \alpha(n)]^{n-1}. \end{aligned}$$ We can then choose $c_{12} \geq 1,$ depending on $m,M,n,$ so that the first inequality of HS4.1(1) holds.  

\medskip

HS4.1(2) follows in \cite{HS79} by the $T$-subharmonicity of the function $\max \{ X_{n+1}-\kappa_{T},0\}^{2},$ together with the fact that $\max \{ X_{n+1} - \kappa_{T},0 \}^{2}|(\B_{2} \cap \spt \partial T) = 0.$ Since both facts hold here as well, then the first inequality of HS4.1(2) holds (with in fact the same $c_{14}$).

\medskip

The proof of the second inequality of R4.1(2) proceeds similarly with only small changes in constants. We take $$c_{15} = 16 \cdot 3^{3n+3} [3+c_{4}+ (M-m) \ALPHA(n-1)][1+M \alpha(n)] c_{16}.$$ We presently may assume \begin{equation} \tag*{R4.1(4)} \E_{\C}(T,1) + \kappa_{T} < 3^{n+2} [1+M \ALPHA(n)] c_{15}^{-1} \sigma^{n+1}, \end{equation} otherwise HS4.1(2) follows directly from R1.6. Using HS2.2(5), HS1.4(1), R1.6(1) as in the proof of HS4.1(5) here gives \begin{equation} \tag*{R4.1(5)} \begin{aligned} \int_{\B_{1}} X_{n+1}^{2} \ d \|T\| & \leq 2 \E_{S}(T,1) + 2 c_{4} \kappa_{T} + 8 \E_{\C}(T,1) \\ & \leq 2 (2+c_{4}+(M-m)\ALPHA(n-1))(E_{\C}(T,1) + \kappa_{T})  \end{aligned} \end{equation} (assuming of course $c_{4} \geq 3$). As we conclude $$\int_{\B_{1}} X_{n+1}^{2} \ d \|T \| \leq c_{15} \sigma^{-n-1} [\E_{\C}(T,1)+\kappa_{T}],$$ then as in \cite{HS79} we must show HS4.1(6) holds. However, using HS4.1(3), Cauchy's inequality, R4.1(5), and R4.1(4) we have $$\begin{aligned} \sup_{\B_{1-\sigma/6} \cap \spt T} X_{n+1}^{2} & \leq 2 \cdot 6^{n} c_{16} \sigma^{-n} \int_{\B_{1}} X_{n+1}^{2} \ d \| T \| + 4 \kappa_{T}^{2} \\ & \leq 4 \cdot 6^{n} [3+c_{4}+(M-m) \ALPHA(n-1)] c_{16} \sigma^{-n} [\E_{\C}(T,1)+\kappa_{T}] \\ & \leq \sigma/12. \end{aligned}$$ The argument that HS4.1(6) holds then proceeds exactly the same.

\medskip

HS4.2. {\it Remark.} Since R4.1(2) remains unchanged, except for that $c_{14},c_{15}$ depend on $m,M,n$, then we conclude from HS2.3(1)(2) that HS4.2(1) holds whenever $|\omega| \leq 1/8,$ $T \in \TT,$ and $$\E_{\C}(T,1)+\kappa_{T} \leq \min \{ (1+c_{4})^{-1},c_{15}^{-1}(1+c_{14})^{-1} 4^{-2n-4} \}.$$ We also have that HS4.2(2)(3)(4) hold, where we still take $$c_{16} = 4^{2n+5}(1+c_{4})(1+c_{12})(1+c_{13})(1+c_{14})(1+c_{15}),$$ which now depends on $m,M,n.$

\medskip

\begin{center}{\bf HS5. Interior nonparametric estimates}\end{center}

\medskip

HS5.1 proves a general decomposition theorem, while HS5.2 and HS5.3 state the well-known gradient estimates for solutions to the minimal surface equation. HS5.4, which proves an approximate graphical decomposition for $T \in \TT$ with sufficiently small cylindrical excess, passes with no serious changes. We introduce in this section $c_{17},\ldots,c_{25}.$

\medskip

HS5.1. LEMMA. We use the exact result here.

\medskip

HS5.2. {\it Remark.} This section introduced standard $L^{2}$ gradient estimates and De Giorgi-Nash H\"older continuity estimates for uniformly elliptic PDEs. Hence, $c_{17},c_{18}$ remain unchanged.

\medskip

HS5.3. {\it Remark.} This section introduces the well-known gradient estimates for solutions to the minimal surface equation, and hence applies HS5.2. We leave $c_{19},c_{20},c_{21},c_{22},c_{23}$ unchanged.

\medskip

HS5.4. THEOREM. The statement passes with no serious change; naturally, we conclude instead the existence of functions $v_{1}^{T} \leq v_{2}^{T} \leq \ldots \leq v_{M}^{T}$ defined over $\V_{T} =\V_{\sigma_{T}} = \V \cap \{ y : \dist (y,\Bdry \V) > \sigma_{T} \}$ and $w_{1}^{T} \leq \ldots \leq w_{m}^{T}$ defined over $\W_{T} = \W_{\sigma_{T}} = \W \cap \{ y: \dist (y,\Bdry \W) > \sigma_{T} \},$ where $\sigma_{T} = c_{24} [\E_{\C}(T,1)+\kappa_{T}]^{1/(2n+3)}.$ Presently $c_{24} \geq 1$ depends on $m,M,n.$

\medskip

Furthermore, HS5.4(1)(2) hold, albeit with $c_{25} \geq 1$ depending on $m,M,n.$ On the other hand, we conclude here, owing to R1.6(1), \begin{equation} \tag*{R5.4(3)} \int_{\V_{T}} \left( \frac{\partial}{\partial r} \left[ \frac{v_{i}^{T}(y)}{|y|} \right] \right)^{2} |y|^{2-n} \ d \LL^{n}y + \int_{\W_{T}} \left( \frac{\partial}{\partial r} \left[ \frac{w_{j}^{T}(y)}{|y|} \right] \right)^{2} |y|^{2-n} \ d \LL^{n}y \end{equation} $$\leq 4 [\E_{S}(T,1)+c_{4} \kappa_{T}] \leq 4 \E_{\C}(T,1) + [(M-m) \ALPHA(n-1) + c_{4} ] \kappa_{T}$$ where $\partial/\partial r [f(y)] = (y/|y|) \cdot Df(y)$ and $c_{4}$ is as in HS2.2.

\medskip

{\it Proof.} The graphical decomposition over $\V_{T},\W_{T}$ depends on 5.3.15 of \cite{F69} and HS5.1, both of which apply in this case. We define $c_{24}$ as in \cite{HS79}, depending on $n,$ $c_{14},$ $c_{15},$ $c_{19},$ $c_{20},$ and thus depending on $m,M,n.$

\medskip

Before proceeding, we remark on a slight technicality in proving HS5.1(1)(2). Applying HS4.1(2) with $\sigma = \sigma_{T}$ gives $$\begin{aligned} \sup_{\V_{T}} |v_{i}^{T}| & \leq \sqrt{2c_{14}c_{15}} \cdot \sigma_{T}^{-n-\frac{1}{2}} [\E_{\C}(T,1)+\kappa_{T}]^{1/2} \\ & \leq \sqrt{2 c_{14} c_{15}} \cdot c_{24}^{-n-\frac{1}{2}} [\E_{\C}(T,1)+\kappa_{T}]^{\frac{1}{2n+3}}. \\ \end{aligned}$$ If $x \in \partial \V_{T},$ then we use HS5.3(1) with $\delta = \dist (x,\Bdry \V_{2 \sigma_{T}/3}) = \sigma_{T}/3,$ giving $$|Dv_{i}^{T}| < \left( \frac{c_{19} \sqrt{2c_{14}c_{15}} [\E_{\C}(T,1)+\kappa_{T}]^{\frac{1}{2n+3}}}{c_{24}^{n+\frac{1}{2}} (\sigma_{T}/3)} \right) \exp \left[\frac{c_{20}\sqrt{2c_{14}c_{15}} [ \E_{\C}(T,1)+\kappa_{T}]^{\frac{1}{2n+3}}}{c_{24}^{n+\frac{1}{2}} (\sigma_{T}/3)} \right] .$$ which is not enough to conclude HS5.4(1). However, HS5.1(1)(2) hold if we instead define $$\V_{T} = \U^{n}(0,1/2) \cap \V \cap \{ y: \dist(y,\Bdry \V) > \sigma_{T} \}$$ $$\W_{T} = \U^{n}(0,1/2) \cap \W \cap \{ y: \dist(y,\Bdry \W) > \sigma_{T} \}$$ In this case we can use HS4.1(2) with $\sigma = 1/2$ to bound $\sup_{\V_{T}} |v_{i}^{T}|$ proportionally to $[\E_{\C}(T,1)+\kappa_{T}]^{1/2}.$ We can thus conclude HS5.4(1)(2) with $c_{25}$ depending on $c_{14}, c_{15}, c_{19}, c_{20},$ and hence on $m,M,n.$ Given what we eventually with to show (in Theorem \ref{main}), this is not a serious issue. 

\medskip

To prove R5.4(3), since HS4.1(2) holds unchanged, then we can verify HS5.1(4)(5) hold exactly. The only difference then comes from the bound given in R1.6(1).

\medskip

\begin{center}{\bf HS6. Blowup sequences and harmonic blowups}\end{center}

\medskip

This section introduces blowup sequences and harmonic blowups, with the aim to prove the necessary rigidity result in HS6.4. Only minor, mostly notational changes must be made. The only serious change is seen in justifying HS6.4(13), which in \cite{HS79} follows from HS3.2, whereas we must use R3.2. We introduce, in HS6.4, $c_{26},\ldots,c_{33}.$

\medskip

HS6.1. We give the same definition of a blowup sequence and harmonic blowup. In this case, we must take functions $v_{i}^{(\nu)},f_{i}$ and $w_{j}^{(\nu)},g_{j}$ respectively with $i \in \{1,\ldots,M\}$ and $j \in \{1,\ldots,m\},$ but still require HS6.1(1)(2)(3)(4) to hold.  

\medskip

As HS4.1(2) holds with no change (except for $c_{14},c_{15}$ now depending on $m,M,n$) then HS6.1(5) also holds here. We also conclude {\it every sequence $S_{1},S_{2},S_{3},\ldots$ in $\TT$ for which} $$\lim_{\nu \rightarrow \infty} [\E_{\C}(S_{\nu},1)+\E_{\C}(S_{\nu},1)^{-1} \kappa_{S_{\nu}}] = 0$$ {\it contains a blowup subsequence.}

\medskip

HS6.2. LEMMA. We conclude the same result, bearing only in mind that in this case we consider functions $f_{1},\ldots,f_{M}$ and $g_{1},\ldots,g_{m}.$

\medskip

{\it Proof.} We take as in \cite{HS79} $$\varepsilon_{\nu} = \E_{\C}(T_{\nu},1)^{1/2}, \ \kappa_{\nu} = \kappa_{T_{\nu}}$$ $$q: \R \times \R^{n+1} \rightarrow \R^{n+1}$$ and define $$\begin{aligned} q(t,(x_{1},\ldots,x_{n+1})) & = (x_{1},\ldots,x_{n-1},t x_{n},t x_{n+1}), \\ Q_{\nu} & = q_{\#} ([0,1] \times [(\partial T_{\nu} \res \C_{2})]) \res \C_{1}. \end{aligned}$$ In this case we take $$S_{\nu} = [T_{\nu}- Q_{\nu}+ (M-m)(\E^{n} \res \W) \times \DELTA_{0} ] \res \C_{1}.$$ By \eqref{appendixA5} and \eqref{appendixA4} with $r=1$ $$\begin{aligned} (\text{Int} \C_{1}) \cap \spt \partial S_{\nu} & = \emptyset, \ \ \p_{\#} S_{\nu} = M \E^{n} \res \U^{n}(0,1), \\ \E_{\C}(S_{\nu},1) & \leq \varepsilon_{\nu}^{2} + \M (Q_{\nu}) \\ & \leq \varepsilon_{\nu}^{2} + \left( \frac{\alpha}{2} \right) (M-m) \ALPHA(n-1) \kappa_{T} \left(1+\frac{\alpha^{2} \kappa_{T}^{2}}{4}+\frac{\alpha^{4} \kappa_{T}^{4}}{16} \right)^{\frac{1}{2}}. \end{aligned}$$ Choosing $S_{\nu,1},\ldots,S_{\nu,M}$ as in HS5.1, then the argument proceeds as in \cite{HS79}.

\medskip

HS6.3. LEMMA. The result holds the same, naturally with $f_{M}$ in place of $f_{m}$ and $g_{m}$ in place of $g_{m-1}.$ Observe as well the typo, we should have instead $U=\U^{n}(a,\sigma).$

\medskip

{\it Proof.} The proof in \cite{HS79} relies on HS4.1(2) to define the boundary data $b_{\nu},$ which is used to solve the minimal surface equation to find the barrier function $u_{\nu}.$ Since HS4.1(2) holds here unchanged, then the proof is exactly the same. 

\medskip

HS6.4. LEMMA. Holds the same, where we must naturally consider $\beta_{1} \leq \beta_{2} \leq \ldots \leq \beta_{M}$ and $\gamma_{1} \geq \gamma_{2} \geq \ldots \geq \gamma_{m}.$ 

\medskip

{\it Proof.} There are some notable changes, particularly in the use of R3.2. 

\medskip

With $\zeta = \max \{|\beta_{1}|,|\beta_{M}|,|\gamma_{1}|,|\gamma_{m}| \}$ and $\delta$ defined the same, since HS5.4(1)(2), HS6.1(1)(2)(3)(4), and HS6.3 hold we can conclude there is $N_{\sigma}$ so that for $\nu \geq N_{\sigma},$ we have that HS6.4(1)(2)(3)(4) hold (where we take $i \in \{ 1,\ldots,M\}$ and $j  \in \{ 1,\ldots,m\}$ respectively in HS6.4(2)(3)). 

\medskip

Next, if $\sigma_{T_{\nu}} < \sigma/4$ then HS5.4(1) implies $|D v_{i}^{(\nu)}| \leq c_{25}$ for $y \in \V_{T_{\nu}}.$ This means we can use HS5.2(3) and then HS5.2(2) as in \cite{HS79} in order to conclude HS6.4(5) for each $i \in \{ 1,\ldots,M \}$ with $c_{26},c_{27}$ depending on $m,M,n.$ We can similarly verify $$\sup_{\W_{\sigma}} |D(w_{j}^{(\nu)}-\varepsilon_{\nu} \gamma_{j} Y_{n})|^{2} \leq c_{27} \ALPHA(n) \sigma^{2} \varepsilon_{\nu}^{2}$$ for each $j \in \{1,\ldots,m\}.$

\medskip

Define $H^{\sigma},I^{\sigma}_{i},J^{\sigma}_{j}$ for $i \in \{1,\ldots,M \},j \in \{1,\ldots,m\},$ $$G^{\sigma}_{\nu} = H^{\sigma} \cup \BETA_{\varepsilon_{\nu}} \left( \bigcup_{i=1}^{M} I_{i}^{\sigma} \right) \cup \BETA_{\varepsilon_{\nu}} \left( \bigcup_{j=1}^{m} J_{j}^{\sigma} \right),$$ and $\Delta^{\sigma}_{\nu}$ as in \cite{HS79}. If we take $\mu : \B^{n}(0,1) \rightarrow [0,1]$ the same $C^{1}$ cut-off function, and define $F_{\nu}^{\sigma}: G_{\nu}^{\sigma} \cup (\R^{n+1} \sim \C_{3/4}) \rightarrow \R^{n+1}$ as in \cite{HS79}, then we again wish to estimate $\M (F^{\sigma}_{\nu \#} T_{\nu}) - \M (T_{\nu}).$ 

\medskip

For this, note the identities HS6.4(6)(7) still hold with $$u_{i}^{(\nu)} = (1-\mu) \varepsilon_{\nu} \beta_{i}(Y_{n} - \sigma) + \mu v_{i}^{(\nu)}.$$ Using HS6.4(2)(5) we get HS6.4(8)(9), with $c_{28},c_{29}$ depending on $m,M,n.$ This implies HS6.4(10)(11), although with $c_{30}$ now depending on $m,M,n.$ 

\medskip

Next, to justify HS6.4(12) we compute using HS1.4(1) $$\begin{aligned} \| T_{\nu} \|(H^{2 \sigma} \cap \C_{3/4+\sigma}) & \leq \M \p_{\#}[T_{\nu} \res H^{2 \sigma} \cap \C_{3/4+\sigma}] + (3/4+\sigma)^{n} \E_{\C}(T,3/4+\sigma) \\ & \leq \M \p_{\#}[T_{\nu} \res H^{2 \sigma} \cap \C_{3/4+\sigma}] + \varepsilon_{\nu}^{2}.\end{aligned}$$ Using \eqref{appendixA5} and \eqref{appendixA4} with $r = 3/4+\sigma$ we can compute $$\begin{aligned} \M \p_{\#}[T_{\nu} \res H^{2 \sigma} \cap \C_{3/4+\sigma}] \leq & 2 \sigma (M+m) \ALPHA(n-1)(3/4+\sigma)^{n-1} \\ & + (M-m) \kappa_{T_{\nu}} \ALPHA(n-1)(3/4+\sigma)^{n+\alpha} \end{aligned}$$ The previous two calculations imply that for sufficiently large $\nu$ \begin{equation} \tag*{R6.4(12)} \| T_{\nu} \| (H^{2 \sigma} \cap \C_{3/4+\sigma}) \leq c_{31} (\sigma + \kappa_{T_{\nu}} + \varepsilon_{\nu}^{2}) \leq 2c_{31} \sigma. \end{equation} for $c_{31}$ depending on $m,M,n.$

\medskip

Combining HS6.4(6)(12) together with R3.2 taking $A = H^{\sigma} \cap \C_{3/4}$ and $\tau = \sigma,$ gives for $c_{32},c_{33}$ depending on $m,M,n$ $$\begin{aligned} \text{R6.4(13)} \ \M F^{\sigma}_{\nu \#} (T_{\nu} \res H^{\sigma}) - \M (T_{\nu} \res H^{\sigma}) & = \M F^{\sigma}_{\nu \#} (T_{\nu} \res H^{\sigma}) - \M (T_{\nu} \res H^{\sigma}) \\ & \leq c_{32} \sigma^{-2} \left[ \kappa_{T_{\nu}} + \int_{H^{2 \sigma} \cap \C_{3/4+\sigma}} X^{2}_{n+1} \ d \| T_{\nu} \| \right] \\ & \leq c_{33} (1+\zeta) \sigma \varepsilon_{\nu}^{2}, \end{aligned}$$ using as well $\kappa_{T_{\nu}} < \sigma^{3} \varepsilon_{\nu}^{2}$ from HS6.4(1). From this we get the same estimate HS6.4(14) for all $\nu \geq N_{\sigma}.$

\medskip

As in \cite{HS79}, the goal is to show the function $\eta : D \rightarrow \R$ is harmonic, where $$\begin{aligned} 
D & = \B^{n}(0,1/2) & \\ \eta (y) & = \beta_{M} Y_{n}(y) &\text{ for } y  \in D \cap \Clos \V \\ \eta (y) & = \gamma_{m} Y_{n}(y) & \text{ for } y  \in D \cap \Clos \W \end{aligned}$$ We proceed in the same way, picking $\sigma_{1},\sigma_{2} \ldots$ any decreasing sequence of numbers with limit zero and $\sigma_{1} < \min \{ \delta/2,1/16 \},$ taking $\nu_{k} = N_{\sigma_{k}}$ and defining $\eta_{k} : D \rightarrow \R$ $$\begin{aligned} 
\eta_{k}(y) & = \beta_{M} (Y_{n}-\sigma_{k})(y) & \text{ for } y  \in D \cap \Clos \V \\ \eta_{k}(y) & = \gamma_{m} (Y_{n}+\sigma_{k})(y) & \text{ for } y  \in D \cap \Clos \W \\ \eta_{k}(y)  & = 0 & \text{ for } y \in D \sim (\V_{\sigma_{k}} \cup \W_{\sigma_{k}}) \end{aligned}$$ Also, define $$\begin{aligned} 
R_{k} & = (-1)^{n-1} [\partial(\E^{n+1} \res \C_{1/2} \cap \{ x: X_{n+1}(x) > \varepsilon_{\nu_{k}} \eta_{k}[\p(x)] \})] \res \C_{1/2} \\ S_{k} & = (-1)^{n-1} [\partial(\E^{n+1} \res \C_{1/2} \cap \{ x: X_{n+1}(x) > \varepsilon_{\nu_{k}} \theta_{k}[\p(x)] \})] \res \C_{1/2}, \end{aligned}$$ (note the slight difference from claimed in \cite{HS79}). We can hence justify HS6.4(15) using HS6.4(7) with $\sigma = \sigma_{k}, \nu = \nu_{k}$ (and the analogous identity over $\W_{\sigma_{k}}$), and HS6.4(1) together with \eqref{appendixA6} with $T = T_{\nu_{k}},$ $\sigma = \sigma_{k}.$ 

\medskip

Finally, HS6.4(16) holds with the same reasoning, with $P_{k} = Q_{\nu_{k}} - F_{\nu_{k} \#}^{\sigma_{k}} Q_{\nu_{k}}$ and $Q_{\nu_{k}} = q_{\#}([0,1] \times [(\partial T_{\nu_{k}}) \res \C_{2}] \res \C_{1})$ as in HS6.2. We also have HS6.4(17), again where $c_{30},c_{33}$ here depend on $m,M,n,$ via HS6.4(14)(15)(16). Therefore, using \eqref{appendixA4} with $r=1,$ we can argue exactly as follows in \cite{HS79} to conclude $\beta_{M} = \gamma_{m},$ as well as the rest of HS6.4.  

\medskip

\begin{center}{\bf HS7. Comparison of spherical and cylindrical excess}\end{center}

\medskip

HS7.1 and HS7.3 give bounds for the cylindrical excess (at smaller radii) in terms of the spherical excess, for $T \in \TT$ with small cylindrical excess. We restate HS7.1 due to a typo in \cite{HS79}. Both hold here with no change. HS7.2 gives a general lemma about homogeneous degree one harmonic functions over $\V.$ We introduce $c_{34},\ldots,c_{37}.$

\medskip

HS7.1. LEMMA. {\it There exist positive constants $c_{34} \geq 1+c_{14},$ $c_{35},$ and $c_{36},$ all depending on $m,M,n,$ so that if $T \in \TT,$} $$\E_{\C}(T,1) + \kappa_{T} \leq c_{34}^{-1},$$ $$\sup_{\C_{1/4} \cap \spt T} X_{n+1}^{2} \leq c_{35}^{-1} \E_{\bS}(T,1),$$ {\it then} $$\E_{\C}(T,1/3) \leq c_{36} [\E_{\bS}(T,1)+\kappa_{T}].$$

\medskip

{\it Proof.} The calculations carry over exactly, although we make a clarification for the reader. 

\medskip

With $c_{4},c_{12},c_{14},c_{15}$ as in HS2.2(5) and HS4.1 (now depending on $m,M,n$), in this case we let $$\begin{aligned} c_{34} & = 2^{2n+2} (1+c_{4})(1+c_{14})(1+c_{15}) \\ c_{35} & = 3^{2n+8} [1+M \ALPHA(n)] c_{12}, \\ c_{36} & = 4^{3n+6} [1+M \ALPHA(n) + c_{4}] c_{12}, \end{aligned}$$ We now assume for contradiction that $T \in \TT$ satisfies the hypothesis of HS7.1, but that $$\E_{\bS}(T,1) + \kappa_{T} < c_{36}^{-1} \E_{\C}(T,1/3).$$ Using HS4.1(2) and $\E_{\C}(T,1)+\kappa_{T} \leq c_{34}^{-1}$ we conclude HS7.1(1).

\medskip

Let $\kappa, \varepsilon, A, \beta, \lambda,$ and $\mu_{k}$ be as in \cite{HS79}. The assumptions $\sup_{\C_{1/4} \cap \spt T} X_{n+1}^{2} \leq c_{35}^{-1} \E_{\bS}(T,1)$ and $\E_{\bS}(T,1)+\kappa_{T} \leq c_{36}^{-1} \E_{\C}(T,1/3)$ then imply $$\sup_{\C_{1/4} \sim \C_{1/8}} |x|^{-1} X_{n+1} \leq 8 c_{35}^{-1/2} \E_{\bS}(T,1)^{1/2} \leq 2c_{36}^{-1/2} \beta \varepsilon.$$ With this fact, we can therefore choose $h_{k}$ a $C^{1}$ vectorfield on $\R^{n+1}$ so that $$\begin{aligned} h_{k}(x) & = 0 & \text{ for } x \in \B_{1/4} \cap \spt T \\ h_{k}(x) & = \lambda^{2}(x) \mu_{k}(|x|) x & \text{ for } x \in \R^{n+1} \sim \B_{1/4}. \end{aligned}$$ As in \cite{HS79}, $h_{k}$ vanishes on $(\spt \partial T) \cup (\B_{1/4} \cap \spt T) \cup \{ x : X_{n+1} \leq |x| (\beta \varepsilon + \kappa) \}.$ The same calculations therefore give $$\int_{A} X_{n+1}^{2} \ d \|T \| \leq 4 (\beta^{2} \varepsilon^{2} + \kappa^{2}) \| T \| (A) + 4^{2n+4} \int_{A} (|X|^{-1} X \cdot \NU^{T})^{2} \ d \|T\|.$$ We also have, using $\sup_{\C_{1/4} \cap \spt T} X_{n+1}^{2} \leq c_{35}^{-1} \E_{\bS}(T,1)$ and $\E_{\bS}(T,1)+\kappa_{T} \leq c_{36}^{-1} \E_{\C}(T,1/3),$ $$\int_{\B_{1/4}} X_{n+1}^{2} \ d \|T\| \leq (c_{35}^{-1} c_{36}^{-1} \varepsilon^{2}) \| T \|(\B_{1/4}) \leq \beta^{2} \varepsilon^{2} \|T\|(\B_{1/4}).$$ Using HS2.3(1), HS4.1(1) with $T,\sigma$ replaced by $(\mu_{3 \#} T) \res \U_{3},1/2,$ HS7.1(1), and HS2.2(5), we conclude $$\begin{aligned} \varepsilon^{2} & \leq 3^{n+2} c_{12} \left[ \kappa + \int_{\C_{1/2}} X_{n+1}^{2} \ d \|T\| \right] \\ & \leq 3^{n+2} c_{12} \left[ \kappa + \int_{\U_{1}} X_{n+1}^{2} \ d \|T\| \right] \\ & \leq 3^{n+2} c_{12} \Big[ \kappa + \beta^{2} \varepsilon^{2} \| T \|(\B_{1/4}) + 4 \beta^{2} \varepsilon^{2} \|T\|(A) \\ & \ \ + 4 \kappa^{2} \|T\|(A) + 4^{2n+4} \int_{A} (|X|^{-1}X \cdot \NU^{T})^{2} \ d \|T\| \Big] \\ & \leq 3^{2n+4} c_{12} [1+ M \ALPHA(n)](\beta^{2} \varepsilon^{2} + \kappa) \\ & \ \ + 4^{3n+6} c_{12} [\E_{S}(T,1) + c_{4} \kappa] \\ & \leq (\varepsilon^{2}/2) + (c_{36}/2) [\E_{S}(T,1) + \kappa] \end{aligned}$$ which is a contradiction.

\medskip

HS7.2. {\it Remark.} HS7.2(1)(2) are general conclusions about homogeneous degree one harmonic functions over $\V,$ which we use exactly.

\medskip

HS7.3. THEOREM. We conclude the same, although with $c_{37}$ depending on $m,M,n.$

\medskip

{\it Proof.} There are only minor changes in the proof.

\medskip

Take for contradiction a sequence $T_{1},T_{2},\ldots$ in $\TT$ satisfying HS7.3(1)(2). Letting $S_{\nu} = (\MU_{3 \#} T_{\nu}) \res \U_{3},$ then HS2.3 implies $S_{\nu} \in \TT.$  With $$\varepsilon_{\nu} = \E_{\C}(S_{\nu},1)^{1/2} = \E_{\C}(T_{\nu},1/3)^{1/2}, \hspace{.25in} \kappa_{\nu} = \kappa_{S_{\nu}} \leq \kappa_{T_{\nu}}/3$$ we can by HS6.1 assume $S_{1},S_{2},\ldots$ is a blowup sequence with associated harmonic blowups $f_{i},g_{j}.$ We can compute using HS1.4(1) and HS2.2(1) $$\begin{aligned} & \limsup_{\nu \rightarrow \infty} \varepsilon_{\nu}^{-2}[\E_{S}(S_{\nu},1)+c_{4} \kappa_{\nu}] \\ & \ \ \leq \limsup_{\nu \rightarrow \infty} (4/3)^{n} \E_{\C}(T_{\nu},1/4)^{-1}[ \exp (c_{1} \kappa_{\nu}) \E_{S}(T_{\nu},1) \\ & \ \ \ \ + (\exp (c_{1} \kappa_{\nu})-1 )\left(\frac{M+m}{2}\right) \ALPHA(n) + c_{4} \kappa_{\nu}] = 0. \end{aligned}$$ Thus, R5.4(3) implies that for all $0 < \rho < 1$ $$f_{i}(\rho y) = \rho f_{i}(y) \text{ for } y \in \V, \hspace{.25in} g_{j}(\rho y) = \rho g_{j}(y) \text{ for } y \in \W.$$

\medskip

Applying HS7.2(1) to $f_{M}-f_{1}$ and HS6.2 to $\min \{ |f_{1}|,\ldots,|f_{M}| \},$ we conclude each $f_{i}$ has zero trace on $\bL.$ Thus, by HS7.2(2), each $f_{i} = \beta_{i} Y_{n}|\V.$ By HS6.2, $\sum_{j=1}^{m} g_{j}$ has zero trace on $\bL$ because $\sum_{i=1}^{M} f_{i}$ does, which implies as in \cite{HS79} that $m|g_{j}|$ has zero trace on $\bL.$ Thus, by HS7.2(2), each $g_{j} = \gamma_{j} Y_{n}| \W.$ HS6.4 then gives $\beta_{1} = \beta_{2} = \ldots = \beta_{M} = \gamma_{1} = \gamma_{2} = \ldots = \gamma_{m}$ and HS7.3(4). 

\medskip

The rest of the proof, which also relies on HS4.1(1)(2), HS4.2, HS7.1, HS7.3(2)(4), is exactly the same. We thus conclude HS7.4(5)(6), which as in \cite{HS79} contradict HS7.4(1).

\medskip

\begin{center}{\bf HS8. Boundary regularity of harmonic blowups}\end{center}

\medskip

This section proves harmonic blowups are given by $C^{2}$ functions over $\V \cup \bL$ and $\W \cup \bL.$ We only need to make a clarification to the end of the proof of HS8.1. We introduce $c_{38},\ldots,c_{43}.$

\medskip

HS8.1. LEMMA. Holds as well, naturally with $i \in \{ 1,\ldots,M\}$ and $j \in \{ 1,\ldots, m\}.$

\medskip

{\it Proof.} For each $\nu \in \{ 1,2,\ldots\}$ letting $\varepsilon_{\nu} = \E_{\C}(T_{\nu},1)^{1/2}$ and $\kappa_{\nu} = \kappa_{T_{\nu}},$ choose for each $0 \leq \sigma \leq 1/12$ an $\omega (\nu,\sigma) \in [-1/8,1/8]$ satisfying HS8.1(1). Applying HS4.1(1)(2), HS4.2(3), HS4.1(1), and HS6.1(1)(2), we infer HS8.1(2).

\medskip

As in \cite{HS79}, suppose $\E_{\C}(\GAMMA_{\omega(\nu,\sigma) \#} \T_{\nu},\sigma/3) \geq \varepsilon_{\nu}^{2}$ for all sufficiently large $\nu.$ By HS8.1(1)(2), HS4.2(3)(4), and HS6.1(1)(2)(3)(4) we can choose a positive integer $N_{\sigma}$ so that, for all integers $\nu \geq N_{\sigma},$ we have $\kappa_{\nu} \leq \varepsilon_{\nu}^{2},$ $$\begin{aligned} S_{\nu} & = (\GAMMA_{\omega (\nu,\sigma) \#} \NU_{\sigma^{-1}} T_{\nu}) \res \U_{3} \in \TT \\ \E_{\C}(S_{\nu},1) + \kappa_{S_{\nu}} & \leq (1/4\sigma)^{n} \E_{\C}((\MU_{4 \#} \GAMMA_{\omega(\nu,\sigma) \#} T_{\nu}) \res \U_{3},1) + \sigma^{\alpha} \kappa_{\nu} \\ & \leq c_{34}^{-1}/2 \\ \E_{\C}(S_{\nu},1/3)+\E_{\C}(S_{\nu},1/3)^{-1} \kappa_{S_{\nu}} & \leq (3/4\sigma)^{n} \E_{\C}((\MU_{4 \#} \GAMMA_{\omega(\nu,\sigma) \#} T_{\nu}) \res \U_{3},1) \\ & \ \ + \varepsilon_{\nu}^{-2} \sigma^{\alpha} \kappa_{\nu} \leq c_{37}^{-1}, \\ \E_{\C}(S_{\nu},1/4) & \leq 2 \E_{C}(\GAMMA_{\eta \#} S_{\nu},1/4) \text{ whenever } |\eta| \leq 1/8. \end{aligned}$$ These are exactly the assumptions of HS7.3 with $T = S_{\nu}$ for $\nu \geq N_{\sigma}$, and so we conclude by HS2.3(2), R1.6(1), and HS6.1(2) $$\begin{aligned} \E_{\C}(\GAMMA_{\omega(\nu,\sigma) \#}, \sigma/4) & = \E_{\C}(S_{\nu},1/4) \leq c_{37} [\E_{\bS}(S_{\nu},1) + \kappa_{S_{\nu}}] \\ & \leq c_{37} [\E_{\bS}(T_{\nu},\sigma) + \sigma^{\alpha} \kappa_{\nu}] \\ & \leq c_{38} [\E_{\bS}(T_{\nu},1) + \kappa_{\nu}] \\ & \leq c_{39} \varepsilon_{\nu}^{2}. \end{aligned}$$ Here, $c_{38}$ depends on $m,M,n,$ as we used $$\E_{\bS}(T_{\nu},\sigma) \leq \exp (c_{1} \kappa_{\nu}) \E_{\bS}(T_{\nu},1) + (\exp(c_{1} \kappa_{\nu})-1) \left(\frac{M+m}{2}\right) \ALPHA(n)$$ by HS2.2(1). Thus, $c_{39}$ also depends on $m,M,n.$

\medskip

Using HS1.4(1) if $\E_{\C}(\GAMMA_{\omega(\nu,\sigma) \#} \T_{\nu},\sigma/3) < \varepsilon_{\nu}^{2},$ we conclude in all cases that $$\E_{\C}(\GAMMA_{\omega(\nu,\sigma) \#} T_{\nu},\sigma/4) \leq c_{40} \varepsilon_{\nu}^{2}$$ for infinitely many $\nu,$ for $c_{40}$ depending on $m,M,n.$ Since $S_{\nu} \in \TT,$ applying HS2.3(1), HS4.1(2) with $T,\sigma$ replaced by $(\GAMMA_{\omega(\nu,\sigma) \#} \MU_{4/\sigma \#} \T_{\nu}) \res \U_{3},1/20,$ and with HS6.1(2) gives HS8.1(3) with $c_{41}$ now depending on $m,M,n.$

\medskip

We conclude by HS8.1(3) and HS6.1(3)(4) that HS8.1(4) holds whenever $\{i,k\} \subset \{ 1,\ldots,M\}$ and $\{j,l\} \subset \{1,\ldots,m\}.$ 

\medskip

Letting $h$ be as in HS6.2 and $H(y) = h(y)-h(\overline{y})$ for $y \in \U^{n}(0,1),$ where $\overline{y} = (y_{1},\ldots,y_{n-1},-y_{n})$ for $y = (y_{1},\ldots,y_{n-1},y_{n}).$ The weak version of the Schwarz reflection principle implies HS8.1(5) for all $0 < \rho < 1.$ 

\medskip

On the other hand, for $y \in \V$ $$\begin{aligned} H & = \sum_{i=1}^{M} f_{i}(y) - \sum_{j=1}^{m} g_{j}(\overline{y}) \\ & = \sum_{i=1}^{M} f_{i}(y) - \sum_{j=1}^{m} (f_{j}(y)+ g_{j}(\overline{y})) + \sum_{j=1}^{m} f_{j}(y) \\ & = 2 \sum_{i=1}^{m} f_{i}(y) + \sum_{i=m+1}^{M} f_{i} - \sum_{j=1}^{m} (f_{j}(y)+ g_{j}(\overline{y})) \\ & = (M+m) f_{k}(y) + 2 \sum_{i=1}^{m} (f_{i}(y)-f_{k}(y)) + \sum_{i=m+1}^{M} (f_{i}(y)-f_{k}(y)) - \sum_{j=1}^{m} (f_{j}(y)+ g_{j}(\overline{y})) \end{aligned}$$ for each $k \in \{1,\ldots,M\},$ and for any $j \in \{1,\ldots,m\}$ $$g_{i}(\overline{y}) = (f_{1}(y)+g_{i}(\overline{y})) - f_{1}(y).$$ We conclude HS8.1 by HS6.1(5) and HS8.1(4)(5).

\medskip

HS8.2. THEOREM. We conclude the same, with $f|\V = f_{1} = f_{2} = \ldots = f_{M}$ and $g|\W = g_{1} = g_{2} = \ldots = g_{m}.$

\medskip

{\it Proof.} The proof is precisely the same, with only minor notational changes.

\medskip

Take a sequence $4 \leq \rho_{k} \rightarrow \infty,$ $f_{i}^{(\rho_{k})}(y) = \rho_{k} f_{i}(y/\rho_{k})$ for $y \in \V,$ $g_{j}^{(\rho_{k})}(y) = \rho_{k} g_{j}(y/\rho_{k})$ for $y \in \W,$ and corresponding $f_{i}^{\ast},g_{j}^{\ast}$ as in \cite{HS79}. By HS5.4(3), HS6.1(2) and HS7.2, each $f_{i}^{\ast},g_{j}^{\ast}$ is some multiple of the function $Y_{n}$ respectively on $\V,\W.$ We wish to show these multiples to coincide. By HS6.4 it is only necessary to show $f_{i}^{\ast},g_{j}^{\ast}$ correspond to a blowup sequence.

\medskip

Define for each $k \in \{1,2,\ldots\}$ $$S_{\nu}^{k} = (\MU_{\rho_{k} \#} T_{\nu}) \res \U_{3},$$ and by HS2.3 and HS6.1(1)(2) an integer $N_{k} \geq k$ so that $S_{\nu}^{k} \in \TT.$ We can as well find $\nu_{k} \geq N_{k}$ so that \begin{equation} \tag*{HS8.2(1)} \begin{aligned} \max \{ \sup_{\V_{1/2}}|f_{1}^{(\rho_{k})}|, \ & \sup_{\V_{1/2}}|f_{M}^{(\rho_{k})}|, \ \sup_{\W_{1/2}}|g_{1}^{(\rho_{k})}|, \ \sup_{\W_{1/2}}|g_{m}^{(\rho_{k})}| \} \\ & \leq \sup_{\C_{1/2} \cap \spt S_{\nu_{k}}^{k}} \varepsilon_{\nu_{k}}^{-1} |X_{n+1}|+k^{-1}, \end{aligned} \end{equation} and so that (applying HS6.3 to $T_{1},T_{2},\ldots$ with $y=0$ and $\sigma = 3/\rho_{k}$) \begin{equation} \tag*{HS8.2(2)} \begin{aligned} & \sup_{\C_{3} \cap \spt S_{\nu_{k}}^{k}} \varepsilon_{\nu_{k}}^{-1} |X_{n+1}| \\ & \leq 3 \max \{ \sup_{\V} |f_{1}^{(\rho_{k}/3)}|, \sup_{\V} |f_{M}^{(\rho_{k}/3)}|, \sup_{\W} |g_{1}^{(\rho_{k}/3)}|,\sup_{\W} |g_{m}^{(\rho_{k}/3)}| \}+k^{-1}, \end{aligned} \end{equation} and so that HS8.2(3) holds with $S_{k}^{\ast} = S_{\nu_{k}}^{k}.$ Using HS4.1(1)(2), we see HS8.2(1)(2) (along with HS2.3(2),HS6.1(2)) imply in case not all $f_{i}^{\ast},g_{j}^{\ast}$ are identically zero that $$0 < \liminf_{k \rightarrow \infty} \varepsilon_{\nu_{k}}^{-1} \E_{\C}(S_{k}^{\ast},1)^{1/2} \leq \limsup_{k \rightarrow \infty} \varepsilon_{\nu_{k}}^{-1} \E_{\C}(S_{k}^{\ast},1)^{1/2} < \infty.$$ This gives there is a positive number $\gamma$ and a blowup subsequence $S_{1}^{\ast},S_{2}^{\ast},\ldots$ whose associated harmonic blowups are $\gamma f_{i}^{\ast},\gamma g_{j}^{\ast}.$

\medskip

Thus, by HS6.2 and HS6.4 we conclude $f_{i}^{\ast} = \beta Y_{n}|\V$ and $g_{j}^{\ast} = \beta Y_{n}|\W$ for some $\beta \in \R.$ Using the Hopf boundary point lemma on $f_{M}-f_{1}$ and $g_{m}-g_{1}$ then concludes the proof.

\medskip

HS8.3. {\it Remark.} The $L^{2}$ estimates given in HS8.3(1)(2)(3) hold with $c_{42}$ the same and $c_{43}$ now depending on $m,M,n.$ We mention that HS8.3(1)(2)(3) follow by the Schwarz reflection principle and well-known estimates for harmonic functions, together with HS6.1(5) taking $\rho = 1/2.$ 

\medskip

\begin{center}{\bf HS9. Excess growth estimate}\end{center}

\medskip

HS9.2 is the central cylindrical excess decay lemma we need to prove HS9.3. Together with the Hopf-type boundary point lemma HS10.1, Theorem \ref{main} follows. All of the results and proofs hold without change. We introduce $c_{44},c_{45},c_{46}.$

\medskip

HS9.1. THEOREM. Holds exactly, but with $c_{44}$ now depending on $m,M,n.$

\medskip

{\it Proof.} The proof is word-for-word the same. The identities used are HS1.4(1), HS6.1, HS8.2, HS8.3(1), HS4.2(4), HS6.3, HS8.3(2)(3), HS4.1(1). 

\medskip

HS9.2. THEOREM. Holds the same.

\medskip

{\it Proof.} The proof is a standard inductive argument using HS9.1, along with HS4.2(3)(4) and HS1.4(1). Follows word-for-word the same.

\medskip

HS9.3. COROLLARY. Holds exactly, except of course we conclude $$\begin{aligned} \p^{-1}(\tilde{V}) \cap \spt \GAMMA_{\eta \#} T & = \bigcup_{i=1}^{M} \graph{\tilde{V}}{\tilde{v}_{i}} \\ \p^{-1}(\tilde{W}) \cap \spt \GAMMA_{\eta \#} T & = \bigcup_{j=1}^{m} \graph{\tilde{W}}{\tilde{w}_{j}} \end{aligned}$$ with $\tilde{V},\tilde{W}$ defined the same and $$\tilde{v}_{1} \leq \tilde{v}_{2} \leq \ldots \leq \tilde{v}_{M}, \hspace{.25in} \tilde{w}_{1} \leq \tilde{w}_{2} \leq \ldots \leq \tilde{w}_{m}.$$

\medskip

{\it Proof.} The proof is virtually the same, so long as we mind to take $i \in \{1,\ldots,M\}$ and $j \in \{1,\ldots,m\}.$ Only HS5.4, HS9.2, and HS2.3(2) are referred to in the proof. Observe that in the process we conclude \begin{equation} \tag*{HS9.3(7)} |D^{2} \tilde{v}_{i}(y)| \leq 2c_{46} |y|^{\beta-1} \end{equation} where $c_{46}$ now depends on $m,M,n,$ with the same holding for $\tilde{w}_{j}.$

\subsection{Concluding Theorem \ref{main}}

Having established HS9.3, the proof of Theorem \ref{main} follows exactly as in the first part of HS11.1.

\medskip

Suppose $T \in \RR^{n}(\U_{3})$ satisfies $(\ast),(\ast \ast)$ from Theorem \ref{main}. Choose $r_{k} \rightarrow 0$ with $r_{k} < 1$ so that as currents $$\MU_{1/r_{k} \#} T \rightarrow M \E^{n} \res \{ (y_{1},\ldots,y_{n}) : y_{n} > 0 \} + m \E^{n} \res \{ (y_{1},\ldots,y_{n}): y_{n} < 0 \}.$$ By 5.4.2 of \cite{F69} $$\lim_{k \rightarrow \infty} \sup_{\B_{r_{k}} \cap \spt T} r_{k}^{-1} X_{n+1} = 0,$$ and so we can choose $k$ sufficiently large so that $$\begin{aligned} T_{k} & = (\MU_{1/r_{k} \#} T) \res \U_{3} \in \TT \\ \max \{ \E_{\C}(T_{k},1), c_{44} \kappa_{T} \} & \leq c_{44}^{-1}. \end{aligned}$$ It follows that HS9.2 holds for $T_{k},$ specifically with $\eta = 0.$ By the Hopf-type boundary point lemma (see Lemma 10.1 of \cite{HS79}) we conclude in applying HS9.2 to $T_{k}$ that $\tilde{v}_{1} = \tilde{v}_{2} = \ldots = \tilde{v}_{M}$ and $\tilde{w}_{1} = \ldots = \tilde{w}_{m}.$ Together with HS9.3(7), we now have Theorem \ref{main}.
  
\appendix

\section{Appendix}

In this section we present some calculations based on the homotopy formula 4.1.9 of \cite{F69}, calculations needed in discussing the sections HS1.6, HS4.1, HS6.2, HS6.4 of \cite{HS79}. We start with the following lemma.

\begin{lemma} \label{appendixlemma1} Suppose $\rho > 0,$ $\sigma<1,$ and that $P \in \RR_{n}(\R^{n})$ nonzero with $$\spt P \cap \{ (y_{1},\ldots,y_{n}) : |(y_{1},\ldots,y_{n-1})| < \rho, |y_{n}| < 1 \} \subset \{ (y_{1},\ldots,y_{n}) : |y_{n}| < \sigma \}$$ and $$\begin{aligned} \partial P \res & \{ (y_{1},\ldots,y_{n}) : |(y_{1},\ldots,y_{n-1})| < \rho, |y_{n}| < 1 \} \\ & = (-1)^{n} \sum_{\ell=1}^{\N} m_{\ell} \Phi_{P,\ell \#}(\E^{n-1} \res \{z : |z| < \rho \}) + (-1)^{n-1} m_{0} \E^{n-1} \res \{ z : |z| < \rho \} \end{aligned}$$ where $m_{\ell}$ are positive integers with $\sum_{\ell=1}^{\N} m_{\ell} = m_{0},$ $$\Phi_{P,\ell}(z_{1},\ldots,z_{n-1}) = (z_{1},\ldots,z_{n-1},\varphi_{P,\ell}(z_{1},\ldots,z_{n-1}))$$ where $\varphi_{P,\ell} \in C^{1}(\R^{n-1} \cap \{z:|z| < 1 \})$ with $\sup_{\R^{n-1} \cap \{z:|z| < 1 \}}|\varphi_{P,\ell}| < \sigma$ for each $\ell = 1,\ldots,\N.$ 

\medskip

Then there are pairwise disjoint open connected sets $O_{P,k}$ and integers $m_{O_{P,k}} \neq 0$ so that $$P \res \{ (y_{1},\ldots,y_{n}) : |(y_{1},\ldots,y_{n-1})| < \rho, |y_{n}| < 1 \} = \sum_{k=1}^{\infty} m_{O_{P,k}} \E^{n} \res O_{P,k}$$ where $m_{O_{P,k}} \in [-m_{0},0)$ if $O_{P,k} \cap \V \neq \emptyset,$ while $m_{O_{P,k}} \in (0,m_{0}]$ if $O_{P,k} \cap \W \neq \emptyset.$ \end{lemma}

\medskip

{\it Proof.} We prove this by induction on $n.$ Note first that the constancy theorem implies \begin{equation} \label{appendixlemma1decomposition} P \res \{ (y_{1},\ldots,y_{n}) : |(y_{1},\ldots,y_{n-1})| < \rho, |y_{n}| < 1 \} = \sum_{k=1}^{\infty} m_{O_{P,k}} \E^{n} \res O_{P,k} \end{equation} for integers $m_{O_{P,k}} \neq 0,$ where $O_{P,k}$ are open connected components of $$\{ (y_{1},\ldots,y_{n}) : |(y_{1},\ldots,y_{n-1})| < \rho, |y_{n}| < \sigma \} \sim \bigcup_{\ell=1}^{\N} \Phi_{P,\ell}(\R^{n-1} \cap \{ z:|z| < \rho \}).$$ We now begin our proof by induction.

\medskip

{\bf n=2:} Define $\tilde{\varphi}_{P,\tilde{\ell}} \in C(\R \cap \{z:|z|<\rho \})$ for $\tilde{\ell} \in \{1,\ldots,\N \}$ with $$\tilde{\varphi}_{P,1} \leq \tilde{\varphi}_{P,2} \leq \ldots \leq \tilde{\varphi}_{P,\N-1} \leq \tilde{\varphi}_{P,\N}$$ and so that for each $z \in \R$ with $|z| < \rho$ we have $$\{ \tilde{\varphi}_{P,\tilde{\ell}}(z) \}_{\tilde{\ell}=1}^{\N} = \{ \varphi_{P,\ell}(z) \}_{\ell=1}^{\N}.$$

\medskip

Take $O_{P,k}$ as in \eqref{appendixlemma1decomposition}, and suppose $O_{P,k} \cap \V \neq \emptyset.$ The constancy theorem, together with $\spt P \cap \{(y_{1},y_{2}):|y_{1}|<\rho, |y_{2}|<1\} \subset \{ (y_{1},y_{2}): |y_{2}| < \sigma \},$ implies there is an open interval $I_{k} \subset (-\rho,\rho)$ and an $\tilde{\ell}_{k} \in \{1,\ldots,\N\}$ so that $$O_{P,k} = \{ (y_{1},y_{2}) : y_{1} \in I_{k}, \ \max \{ 0,\tilde{\varphi}_{P,\tilde{\ell}_{k}-1}(y_{1})\} < y_{2} < \tilde{\varphi}_{P,\tilde{\ell}_{k}}(y_{1}) \}.$$

\medskip

First, suppose $\tilde{\varphi}_{P,\tilde{\ell}_{k}}(y_{1}) = \tilde{\varphi}_{P,\N}(y_{1})$ for each $y_{1} \in I_{k}.$ It follows we can find $\ell_{1},\ldots,\ell_{\N-\tilde{\ell}_{k}+1} \in \{1.\ldots,\N \}$ so that $$\varphi_{P,\ell_{1}}(y_{1}) = \ldots = \varphi_{P,\ell_{N-\tilde{\ell}_{k}+1}}(y_{1}) = \tilde{\varphi}_{P,\N}(y_{1})$$ for each $y_{1} \in I_{k},$ and hence $$\partial O_{P,k} \cap \{ (y_{1},y_{2}): \max \{ 0,\tilde{\varphi}_{P,\tilde{\ell}_{k}-1}(y_{1}) \} < y_{2} \} = \Phi_{P,\ell_{1}}(I_{k}) = \ldots = \Phi_{P,\ell_{\N-\tilde{\ell}_{k}+1}}(I_{k}).$$ From this one can show $m_{O_{P,k}} = - (m_{\ell_{1}}+\ldots+m_{\ell_{N-\tilde{\ell}_{k}+1}}),$ in order to conclude $m_{O_{P,k}} \in [-m_{0},0).$

\medskip

Second, suppose $\tilde{\varphi}_{P,\tilde{\ell}_{k}}(y_{1}) = \tilde{\varphi}_{P,\N-1}(y_{1})$ for each $y_{1} \in I_{k},$ but $\tilde{\varphi}_{P,\N-1}(y_{1}) < \tilde{\varphi}_{P,\N}(y_{1})$ for some $y_{1} \in I_{k}.$ We can thus find an open interval $\tilde{I}_{k} \subset I_{k}$ and an $O_{P,\tilde{k}}$ from \eqref{appendixlemma1decomposition} disjoint from $O_{P,k}$ so that $$\partial O_{P,\tilde{k}} \cap \{ (y_{1},y_{2}) : y_{1} \in \tilde{I}_{k},y_{2}>0 \} = \{ (y_{1},\tilde{\varphi}_{P,\N-1}(y_{1})): y_{1} \in \tilde{I}_{k} \} \cup \{ (y_{1},\tilde{\varphi}_{P,\N}(y_{1})) : y_{1} \in \tilde{I}_{k} \}.$$ Using the previous case applied to $O_{P,\tilde{k}}$ implies there are $\ell_{1},\ldots,\ell_{\N-\tilde{\ell}_{k}+1} \in \{1,\ldots,\N \}$ so that $$m_{O_{P,k}} = - (m_{\ell_{1}}+\ldots+m_{\ell_{\N-\tilde{\ell}_{k}+1}}),$$ and hence again $m_{O_{P,k}} \in [-m_{0},0).$

\medskip

Third, we can argue inductively that every $m_{O_{P,k}} \in [-m_{0},0)$ whenever $O_{P,k} \cap \V \neq \emptyset.$ By likewise first considering $\tilde{\varphi}_{P,1},$ we can show $m_{O_{P,k}} \in (0,m_{0}]$ whenever $O_{P,k} \cap \W \neq \emptyset.$ This shows the case $n=2.$

\medskip

{\bf n $>$ 2.} Take any of almost every $t \in (-\rho,\rho)$ such that the slice $$<P,X_{1},t> = \partial [P \res \{ X_{1} < t \}] - (\partial P) \res \{ X_{1} < t \}.$$ exists, by 4.3.6 of \cite{F69}. Note that $$\begin{aligned} \spt <P,X_{1},t> \cap \{ (y_{1},\ldots,y_{n}) & : |(y_{2},\ldots,y_{n-1})| < \rho-t, |y_{n}| < 1 \} \\ & \subset \{ (t,y_{2},\ldots,y_{n}): |(y_{2},\ldots,y_{n-1})| < \rho-t, |y_{n}| < \sigma \}. \end{aligned}$$ and $$\begin{aligned} \partial <P,X_{1},t> & \res \{ (t,y_{2},\ldots,y_{n}) : |(y_{2},\ldots,y_{n-1})| < \rho-t, |y_{n}| < 1 \} \\ & = \delta_{t} \times \left[ (-1)^{n-1} \sum_{\ell=1}^{\N} m_{\ell} \Phi^{t}_{P,\ell \#}(\E_{n-2}) + (-1)^{n-1} m_{0} \E^{n-2} \right] \end{aligned}$$ where $\Phi_{P,\ell}^{t}(z_{1},\ldots,z_{n-2}) = (z_{1},\ldots,z_{n-2},\varphi_{T,\ell}(t,z_{1},\ldots,z_{n-2})).$ By induction, with $\rho$ replaced by $\rho-t,$ it follows we can write $$<P,X_{1},t> \res \{ (t,y_{2},\ldots,y_{n}) : |(y_{2},\ldots,y_{n-1})| < \rho-t, |y_{n}| < 1 \} = \delta_{t} \times \sum_{k=1}^{\infty} m_{O^{t}_{P,k}} \E^{n-2} \res O^{t}_{P,k}$$ where $O^{t}_{P,k} \subset \R^{n-1}$ are open connected sets so that $m_{O^{t}_{P,k}} \in [-m_{0},0)$ whenever $O^{t}_{P,k} \cap \{ (z_{1},\ldots,z_{n-1}):z_{n-1}>0 \} \neq \emptyset,$ and $m_{O_{t,k}} \in (0,m_{0}]$ whenever $O_{t,k} \cap \{ (z_{1},\ldots,z_{n-1}):z_{n-1}<0 \} \neq \emptyset.$

\medskip

Therefore, since each $m_{O_{P,k}}$ as in \eqref{appendixlemma1decomposition} is $m_{O_{P,k}}= m_{O^{t}_{P,\tilde{k}}}$ for some $t \in (-\rho,\rho),$ then the lemma holds. 

\medskip

\begin{lemma} \label{appendixlemma2} Let $q: \R \times \R^{n+1} \rightarrow \R^{n+1}$ be given by $$q(t,x_{1},\ldots,x_{n+1}) = (x_{1},\ldots,x_{n-1},tx_{n},tx_{n+1}).$$ For $T \in \TT,$ define $$Q_{T} = q_{\#}([0,1] \times (\partial T \res \C_{2})).$$ Then for every $r \in (0,2)$ \begin{equation} \label{appendixA4} \begin{aligned} \M(Q_{T} \res \C_{r}) & \leq \left( \frac{\alpha}{2} \right) (M-m) \kappa_{T} \ALPHA(n-1) \left( 1+\frac{\alpha^{2} \kappa_{T}^{2}}{4} r^{2 \alpha}+ \frac{\alpha^{4} \kappa_{T}^{4}}{16} r^{4 \alpha} \right)^{\frac{1}{2}} r^{n+\alpha} \\ \M \p_{\#} Q_{T} \res \C_{r} & \leq (M-m) \kappa_{T} \ALPHA(n-1) r^{n+\alpha}. \end{aligned} \end{equation} Furthermore, \begin{equation} \label{appendixA5} \p_{\#} T \res \C_{r} = (M \E^{n} \res \V + m \E^{n} \res \W + \p_{\#} Q_{T}) \res \C_{r}. \end{equation} For any $\sigma \in (0,1/2),$ if we have $\kappa_{T} < \sigma,$ then \begin{equation} \label{appendixA6} \begin{aligned} \M ((\p_{\#} T - \E^{n}) \res & \U^{n}(0,1/2) \cap \{ (y_{1},\ldots,y_{n}) : |y_{n}| < \sigma \} ) \\ = & \M ( \p_{\#} T \res \U^{n}(0,1/2) \cap \{ (y_{1},\ldots,y_{n}) : |y_{n}| < \sigma \}) \\ & - \M (\E^{n} \res \U^{n}(0,1/2) \cap \{ (y_{1},\ldots,y_{n}) : |y_{n}| < \sigma \}).\end{aligned} \end{equation} \end{lemma}

\medskip

{\it Proof.} First, we compute by 4.1.9 of \cite{F69} (see also the end of R3.2) $$\begin{aligned} \M(Q_{T} \res \C_{r}) & \leq \| \partial T \res \C_{r} \| (\sqrt{X_{n}^{2}+X_{n+1}^{2}}) \\ & \leq \left( \frac{\alpha}{2} \kappa_{T} r^{1+\alpha} \right) \cdot \| \partial T \| \C_{r} \\ & \leq \left( \frac{\alpha}{2} \right)(M-m) \kappa_{T} \ALPHA(n-1) \left( 1+\frac{\alpha^{2} \kappa_{T}^{2}}{4} r^{2\alpha}+ \frac{\alpha^{4} \kappa_{T}^{4}}{16} r^{4 \alpha} \right)^{\frac{1}{2}} r^{n+\alpha}. \end{aligned}$$ A similar calculation, for $\M \p_{\#} Q_{T} \res \C_{r},$ can be used to conclude \eqref{appendixA4}

\medskip

Second, observe that for any $r \in (0,2)$ we have by 4.1.8-9 of \cite{F69} $$\partial ( \p_{\#} T - \p_{\#} Q_{T} + (M-m) \E^{n} \res \W) \res \C_{r} = 0,$$ since $\sum_{\ell=1}^{\N} m_{\ell} = M-m.$ This proves \eqref{appendixA5}, by the constancy theorem and the definition of $\TT.$

\medskip

Third, if $\kappa_{T} < \sigma$ then we can apply Lemma \ref{appendixlemma1} with $P = (\p_{\#} Q_{T}) \res \U^{n}(0,1),$ $\rho = 1/2,$ and $m_{0} = M-m$ to get $$P \res \{ (y_{1},\ldots,y_{n}): |(y_{1},\ldots,y_{n-1})|<1/2, |y_{n}| < 1 \} = \sum_{k=1}^{\infty} m_{O_{P,k}} \E^{n} \res O_{P,k},$$ where $m_{O_{P,k}} \in [-(M-m),0)$ if $O_{P,k} \cap \V \neq \emptyset,$ while $m_{O_{P,k}} \in (0,(M-m)]$ if $O_{P,k} \cap \W \neq \emptyset.$ By \eqref{appendixA5}, we see that \eqref{appendixA6} follows.

\medskip

\begin{lemma} \label{appendixlemma3} Suppose $T \in \TT,$ and for $\tau \in (0,1)$ suppose $\kappa_{T} \leq 4 \cdot 3^{n} [1+M\ALPHA(n)] \tau^{2}.$ Let $\mu \in C^{1}(\R^{n})$ be a function so that $$\begin{aligned} & \mu(y) = 0 & \text{ for } & |y| \leq 1 \\ 0 < & \mu(y) <1 & \text{ for } 1 \leq & |y| \leq 1 + \tau \\ & \mu(y) = 1 & \text{ for } & |y| > 1+\tau \end{aligned}$$ and with $|D \mu| \leq 3/\tau.$ Define $F : \R^{n+1} \rightarrow \R^{n+1}$ by $$F(x_{1},\ldots,x_{n+1}) = (x_{1},\ldots,x_{n},\mu(x_{1},\ldots,x_{n}) x_{n+1}).$$ If we let $h : \R \times \R^{n+1} \rightarrow \R^{n+1}$ be given by $h(t,x) = (1-t)F(x)+tx$ and $R_{T} = h_{\#} ([0,1] \times \partial T),$ then \begin{equation} \label{appendixA8} \M(R_{T}) \leq \frac{1}{2} \left[ \frac{\sqrt{21}}{4}+ 2^{\frac{9n-5}{2}} 3^{n^{2}-\frac{1}{2}} \right] (M-m) \ALPHA(n-1) [1+M\ALPHA(n)]^{n-1} \kappa_{T}. \end{equation} \end{lemma} 

{\it Proof.} We compute using 4.1.9 of \cite{F69} (see also the end of R3.2), $$\begin{aligned} \M(R_{T}) & \leq \| \partial T \| \left( (1-\mu)|X_{n+1}| \sup \left\{ 1, \left(1+\mu^{2}+ X_{n+1}^{2} |D \mu|^{2} \right)^{\frac{n-1}{2}} \right\} \right) \\ & \leq \left( \frac{\alpha}{2} \kappa_{T} \right) \| \partial T \| (\C_{1}) \\ & \ \ + \left( \frac{\alpha}{2} \kappa_{T} \right) (1+\tau)^{1+\alpha} \left(2+ \left( \frac{\alpha^{2}}{4} \right) \kappa_{T}^{2} (1+\tau)^{2+2\alpha} \left( \frac{9}{\tau^{2}}\right) \right)^{\frac{n-1}{2}} \| \partial T \|(\C_{1+\tau} \sim \C_{1}) \\ & \leq \frac{1}{2} \left[ \| \partial T \|(\C_{1}) + 2^{2} \left(2+ \frac{36 \kappa_{T}^{2}}{\tau^{2}} \right)^{\frac{n-1}{2}} \| \partial T \|(\C_{1+\tau} \sim \C_{1}) \right] \kappa_{T} \\ & \leq \frac{1}{2} \left[ \| \partial T \|(\C_{1}) + 2^{2} \left( 2+ 2^{6} \cdot 3^{2n+2} [1+M \ALPHA(n)]^{2} \right)^{\frac{n-1}{2}} \| \partial T \|(\C_{1+\tau} \sim \C_{1}) \right] \kappa_{T} \\ & \leq \frac{1}{2} \left[ \| \partial T \|(\C_{1}) + 2^{\frac{7n-3}{2}} \cdot 3^{n^{2}-1} [1+M\ALPHA(n)]^{n-1} \| \partial T \|(\C_{1+\tau} \sim \C_{1}) \right] \kappa_{T} \\ & \leq \frac{1}{2} \left[ \frac{\sqrt{21}}{4}+ 2^{\frac{9n-5}{2}} 3^{n^{2}-\frac{1}{2}} \right] (M-m) \ALPHA(n-1) [1+M\ALPHA(n)]^{n-1} \kappa_{T}. \\ \end{aligned}$$

\end{flushleft}
\end{document}